\documentclass[11pt,reqno]{amsart}
\usepackage{amsmath, amsthm, amsfonts, amssymb, color}
\usepackage{mathrsfs, bbm}

\allowdisplaybreaks

\newtheorem{thm}{Theorem}[section]
\newtheorem{cor}[thm]{Corollary}
\newtheorem{lem}[thm]{Lemma}
\newtheorem{remark}[thm]{Remark}
\newtheorem{prp}[thm]{Proposition}

\theoremstyle{definition}
\newtheorem{defn}{Definition}[section]

\definecolor{wco}{rgb}{0.5,0.2,0.3}

\numberwithin{equation}{section} \theoremstyle{remark}

\def\1{{\mathbbm 1}}

\def\wt{\widetilde}

\def\R{\mathbb R}

\def\E{\mathbb E}

\def\L{\mathcal L}

\def\<{\langle}
\def\>{\rangle}

\def\pf{\noindent{\bf Proof.} }

 \def\beq{\begin{equation}}

 \def\P{\mathbb P} 
  \def\ee{\varepsilon}

\begin{document}
\bibliographystyle{plain}

\title{ Dirichlet heat kernel estimates for  isotropic L\'evy processes with Gaussian components in Lipschitz open sets }

\author{Jie-Ming Wang}
\address{School of Mathematics and Statistics,  Beijing Institute of Technology,
Beijing 100081, P. R. China}
\curraddr{}
\email{wangjm@bit.edu.cn}
\thanks{}



\date{}

\begin{abstract}
In this paper,  the two-sided Dirichlet heat kernel estimates are obtained for a class of discontinuous isotropic L\'evy processes with Gaussian components   in Lipschitz open sets.
Furthermore,  the necessary and sufficient conditions for the Varopoulos-type  Dirichlet heat kernel estimates holding for such processes in Lipschitz open sets are derived.
\end{abstract}

\maketitle

\section{Introduction}

The  study of Dirichlet heat kernel and its estimates takes up an important place in both  analysis and probability theory.
 Dirichlet heat kernel for an operator in an open set $D$ is the fundamental solution of the corresponding heat equation for the operator  with zero exterior condition.
When $X$ is a Markov process with the infinitesimal generator $\L,$
the transition density function  of the subprocess of $X$ killed upon leaving $D$ (if it exists) is  the fundamental solution for the operator $\L$ with zero exterior condition, which is also called the Dirichlet heat kernel for $X$ in $D$.
The Dirichlet heat kernel estimates for the Laplace operator in $C^{1, 1}$ domains are established in Davies \cite{D1, D2} and Davies and Simon \cite{DS} for the upper bound estimates and in Zhang \cite{Zh} for the lower bound estimates.
 Varopoulos \cite{V} obtained the qualitative estimates for the Dirichlet heat kernel $q_D(t, x, y)$ for  differential operators  in  Lipschitz open sets, which is shown to be comparable to the global transition density function in $\R^d$ multiplied by the survival probabilities $\P_\cdot(\tau_D>t)$ with starting points $x$ and $y$ for a Lipschitz open set $D,$ where $\tau_D$ is the first exit time of the diffusion from $D.$
The Dirichlet heat kernel estimates for diffusion in inner uniform domains on metric measure spaces are established in \cite{GSC, LS}.

For pure jump processes and the associated non-local operators, Chen, Kim and Song \cite{CKS1} obtained the explicit two-sided Dirichlet heat kernel estimates for the fractional Laplacian $-(-\Delta)^{\alpha/2} (\alpha\in (0, 2))$ in $C^{1, 1}$ open sets in $\R^d$.
After this, many significant progresses have been made on the Dirichlet heat kernel estimates for more pure jump  processes  in $C^{1, 1}$ open sets, including censored stable-like processes, relativistic stable processes, pure jump subordinate Brownian motion, a class of unimodal L\'evy processes as well as a class of   symmetric pure jump processes which are not necessarily L\'evy processes, see \cite{BGR2, CKS2, CKS3, CKS5, CKS6, GKK, KM} etc.. In the case of non-smooth open sets, Bogdan,  Grzywny and Ryznar \cite{BGR1} obtained the Varopoulos-type   Dirichlet heat kernel estimates for the fractional Laplacian  operator $-(-\Delta)^{\alpha/2} (\alpha\in (0, 2))$ in $\kappa$-fat open sets.  The Varopoulos-type  Dirichlet heat kernel estimates for   a   class of purely discontinuous rotationally  symmetric L\'evy  processes in $\kappa$-fat open sets are shown to hold in  Chen, Kim and Song \cite{CKS6}. In  \cite{ChoKSV},  the result is obtained for a class of  purely non-local operators with killing potentials in $\kappa$-fat open sets.

 Although the significant progresses on pure jump processes, there are limited results of Dirichlet heat kernel estimates for discontinuous  L\'evy processes with Gaussian components at present, especially on non-smooth open sets.  Chen, Kim and Song \cite{CKS4}  established the explicit  Dirichlet heat kernel estimates for $\Delta+\Delta^{\alpha/2}(\alpha\in (0, 2))$ in $C^{1, 1}$ open sets in $\R^d.$
This result has been later extended to  subordinate Brownian motions with Gaussian components in $C^{1, 1}$ open sets (see \cite{CKS7, BK})  as well as a class of subordinate diffusions with diffusive components in $C^{1, s} (s\in (0, 1))$ open sets in \cite{W}. In the case of non-smooth open sets, the question whether the Varopoulos-type  Dirichlet heat kernel estimates holds for $\Delta+\Delta^{\alpha/2}(\alpha\in (0, 2))$   in Lipschitz open sets has remained unanswered.

\smallskip

The purpose of this paper is to investigate this question and establish the two-sided Dirichlet heat kernel estimates for a class of isotropic L\'evy processes with Gaussian components in Lipschitz open sets.
More specifically,  consider the operator
\begin{equation}\label{e:L}
\L f(x)=\Delta f(x)+\lim_{\ee\rightarrow 0}\int_{\{z\in \R^d: |z|>\ee\}} (f(x+z)-f(x))j(z)\,dz, \quad f\in  C_c^\infty(\R^d),
\end{equation}
where $j(z)=j(|z|):=\dfrac{1}{|z|^d \phi(|z|)}$ and
 $\phi$ is a strictly increasing continuous
function $ \R_+\mapsto \R_+$ with $\phi (0)=0$ and $\phi(1)=1$
satisfying there are constants $c\ge 1$, $0<\beta_1\leq \beta_2 <2$
such that
\begin{equation}\label{e:1.1}
c^{-1} \Big(\frac Rr\Big)^{\beta_1} \leq   \frac{\phi (R)}{\phi (r)}
\leq c \, \Big(\frac Rr\Big)^{\beta_2} \qquad \hbox{for every }
0<r<R<\infty.
\end{equation}
 Observe that  condition \eqref{e:1.1}  implies that
\begin{equation}\label{e:1.3}
c^{-1} r^{\beta_1} \leq \phi (r) \leq c r^{\beta_2}
\qquad \hbox{for } r\geq 1
\end{equation}
and
\begin{equation}\label{e:1.4}
c^{-1} r^{\beta_2} \leq \phi (r) \leq c r^{\beta_1}
\qquad \hbox{for } r\in (0, 1].
\end{equation}
In particular, $\L=\Delta+\Delta^{\alpha/2}$ when $\phi(z)=\mathcal{A}(d, \alpha)^{-1}|z|^\alpha$ for $\alpha\in (0, 2),$ where $\mathcal{A}(d, \alpha)$ is a positive normalizing constant so that the Fourier transform of $\widehat{\Delta^{\alpha/2}f}$ of $\Delta^{\alpha/2}f$ is $-|\xi|^\alpha \hat f(\xi).$ It follows from  Chen and Kumagai \cite[Theorem 1.2]{CK3} that  there is a discontinuous L\'evy process with Gaussian component $X$ in $\R^d$ associated  with the generator $\L$  and
  $X$ has a jointly continuous transition density function $p(t, x, y)$ with respect to the Lebesgue measure in $\R^d.$

  The following two-sided estimates of $p(t, x, y)$ are established in \cite[Theorem 1.4]{CK3}  or \cite[Theorem 1.4]{CKKW}.

\begin{thm}\label{T1}
There are positive constants $C_k=C_k(d, \phi)$, $k=1, 2, 3$
such that for every $t>0$ and $x, y \in \R^d$,
\begin{eqnarray}
&& C_1^{-1}\,    \left( t^{-d/2}\wedge\phi^{-1}(t)^{-d} \right)   \wedge
\left( p^c(t, C_2|x- y|)+p^j(t, |x- y|) \right)
\nonumber\\
&\leq & p(t, x, y) \leq C_1\, \left( t^{-d/2}\wedge\phi^{-1}(t)^{-d}
\right) \wedge \left( p^c(t, C_3 |x- y|)+p^j(t, |x- y|) \right).
\label{eq:HKj+dul}\end{eqnarray}
where
\begin{equation}\label{eqn:4}
p^c(t, r):=t^{-d/2}e^{-\frac{r^2}{2t}} \qquad \hbox{and} \quad  p^j(t, r) :=\left( \phi^{-1}(t)^{-d} \wedge
\frac{t}{r^d \phi (r)} \right)
\end{equation}
with $\phi^{-1}$ being the inverse function of $\phi$.
\end{thm}

\smallskip

The objective of  this paper is to obtain  the two-sided Dirichlet heat kernel estimates for  $X$  in Lipschitz open sets in $\R^d (d\geq 1).$
The result shows the Varopoulos-type Dirichlet heat kernel estimates  fails to hold even for $\Delta+\Delta^{\alpha/2} (\alpha\in (0, 2))$ in some Lipschitz open sets.  This distinguishes the boundary behavior for the non-local operators with Gaussian components from that for the elliptic operator and the purely non-local operators.
The necessary and sufficient conditions for the Varopoulos-type  Dirichlet heat kernel estimates holding for $X$ in Lipschitz open sets  are obtained in this paper.

\smallskip

\smallskip

The following is the first main theorem of this paper.  For each open set $B,$ define $\tau_B:=\inf\{t>0: X_t\notin B\}$ the first exiting time of $X$ from $B.$ Let $D$ be an  open set in $\R^d.$
  Denote by $p_D(t, x, y)$ the transition density function of the subprocess of $X$ killed upon exiting $D.$

\begin{thm}\label{T2}
Suppose that $D$ is a Lipschitz  open set in $\R^d$ with characteristics $(R_0, \Lambda_0).$

{\rm (i)} For each $T>0,$ there exists a positive constant $C_1=C_1(d, \phi, R_0, \Lambda_0, T)$  such that for   $t\in (0, T)$ and  $x, y\in D$ with $|x-y|\leq t^{1/2},$
$$
C_1^{-1} \P_x(\tau_D>t)t^{-d/2}\P_y(\tau_D>t) \leq p_D(t, x, y)\leq C_1 \P_x(\tau_D>t)t^{-d/2}\P_y(\tau_D>t).
$$

{\rm (ii)}
Assume  the path distance in each connected component of $D$ is comparable to the Euclidean distance with characteristic $\chi_1.$ For each $T>0,$
there exist positive constant $C_2=C_2(d, \phi, R_0, \Lambda_0, \chi_1, T)$ and $C_k=C_k(d, \phi, \chi_1), k=3, 4$ such that for   $t\in (0, T)$ and  $x, y\in D$ with $|x-y|> t^{1/2},$
$$
\begin{aligned}
&C_2^{-1} \P_x(\tau_D>t)p(t, C_3x, C_3y)\P_y(\tau_D>t)+C_2^{-1}\int_0^t \P_x(\tau_D>s)\P_y(\tau_D>t-s)\,ds \cdot j(|x-y|) \\
&\qquad \leq p_D(t, x, y)\\
&\leq C_2 \P_x(\tau_D>t)p(t, C_4x, C_4y)\P_y(\tau_D>t)+C_2\int_0^t \P_x(\tau_D>s)\P_y(\tau_D>t-s)\,ds \cdot j(|x-y|).
\end{aligned}
$$

{\rm (iii)} Suppose in addition that $D$ is bounded. Then there exists a constant $C=C(d, \phi, R_0, \Lambda_0, {\rm diam}(D))$ such that for $t\geq 3$ and $x, y\in D,$
$$C^{-1}e^{-\lambda_1 t}\P_x(\tau_D>1)\P_y(\tau_D>1)\leq p_D(t, x, y)\leq Ce^{-\lambda_1 t}\P_x(\tau_D>1)\P_y(\tau_D>1),$$
where $-\lambda_1<0$ is the largest eigenvalue of the generator of $X^D.$

\end{thm}

See Definitions \ref{D1} and \ref{D2} for the definitions of a Lipschitz open set with characteristics $(R_0, \Lambda_0)$ and the path distance in each connected component of $D$ is comparable to the Euclidean distance with characteristic $\chi_1.$
Theorem \ref{T2} shows that the Varopoulos-type estimate holds for $p_D(t, x, y)$ in the   near diagonal case $|x-y|\leq t^{1/2}.$
 However, this may not hold for $p_D(t, x, y)$ in the off-diagonal case $|x-y|> t^{1/2}.$

For each open set $B,$ denote by  $G_B(x, y)$ the Green function of $X$ killed upon exiting $B.$
Let $D$ be a Lipschitz open set with characteristics $(R_0, \Lambda_0).$ It is  known (see e.g. \cite[Lemma 6.6]{ChZ}) that
there exists   $ \kappa=\kappa(d, \Lambda_0)\in (0, 1/4)$ such that for $r\in (0, R_0)$
and $z\in\partial D,$
\begin{equation}\label{e:1.8}
\hbox{there exists } z_r\in D\cap \partial B(z, r)  \hbox{ with }
\kappa r\leq \delta_D(z_r)< r.
\end{equation}
In the following result, the necessary and sufficient conditions for the Varopoulos-type Dirichlet heat kernel estimates holding for the process $X$ in a Lipschitz open set $D$ will be given.

\begin{thm}\label{T3}
Suppose that $D$ is a Lipschitz open set in $\R^d$ with characteristics $(R_0, \Lambda_0).$
 Assume  the path distance in each connected component of $D$ is comparable to the Euclidean distance with characteristic $\chi_1.$
The followings are equivalent:

\smallskip

{\rm (i)} The Varopoulos-type estimate holds for $p_D(t, x, y)$ for $x, y\in D$ and $t\in (0, T).$
That is, for each $T>0,$ there exist positive constants $C_1=C_1(d, \phi, R_0, \Lambda_0, \chi_1, T)$ and $C_k=C_k(d, \phi,  \chi_1), k=2, 3$ such that  for  any $x, y\in D$ and $t\in (0, T),$
\begin{equation}\label{e:1.9}
C_1^{-1}\P_x(\tau_D>t)p(t, C_2x, C_2y)\P_y(\tau_D>t)\leq p_D(t, x, y)\leq C_1\P_x(\tau_D>t)p(t, C_3x, C_3y)\P_y(\tau_D>t).
\end{equation}

\smallskip

{\rm (ii)} There exists $C=C(d, \phi, R_0, \Lambda_0)>1$ such that  for  any $z_0\in \partial D$ and $t\in (0, R_0^2),$
\begin{equation}\label{e:1.10}
C^{-1} \E_x \tau_{D\cap B(z_0, \sqrt t)}\leq G_{D\cap B(z_0, \sqrt t)}(x, x_t)t^{d/2}\leq C \E_x \tau_{D\cap B(z_0, \sqrt t)} \quad \hbox{for }\: x\in  D\cap B(z_0, \sqrt t/4),
\end{equation}
  where $x_t$ is a point of $D\cap \partial B(z_0, \sqrt t/2)$ with $\delta_D(x_t)\geq \kappa\sqrt t/2$ and $\kappa$ is the constant in \eqref{e:1.8}.

\smallskip

{\rm (iii)} The scale invariant boundary Harnack principle holds for $X$ in Lipschitz open set $D$.
That is,  there exists $ C= C(d, \phi, R_0,  \Lambda_0)>0$  such that for any $z_0\in\partial D, r\in (0, R_0/2)$ and any two  nonnegative harmonic functions $h_1$ and $h_2$ with respect to   $X$  in
 $D\cap B(z_0, 2r)$  vanishing continuously on $D^c\cap B(z_0, 2r),$
\begin{equation}\label{e:1.13}
\dfrac{h_1(x)}{h_1(y)}\leq C \dfrac{h_2(x)}{h_2(y)} \quad \hbox{for any }
x,y\in D\cap B(z_0, r).
\end{equation}
\end{thm}

\begin{remark} \rm

(i) It is known from \cite{CW2} that  the scale invariant boundary Harnack principle  holds for a large class of subordinate Brownian motion with Gaussian components in Lipschitz domains in $\R^d$ satisfying the interior cone condition with common angle
  $\theta\in (\cos^{-1}(1/\sqrt d), \pi)$,
but fails  in any truncated circular cone  with angle $\theta \leq \cos^{-1}(1/\sqrt d)$.
In view of Theorem \ref{T3}, as a particular case, this  indicates that the Varopoulos-type Dirichlet heat kernel estimates \eqref{e:1.9}   holds for $\Delta+\Delta^{\alpha/2}$ in Lipschitz domains in $\R^d$ satisfying the interior cone condition with common angle
  $\theta\in (\cos^{-1}(1/\sqrt d), \pi)$,
but fails  in any truncated circular cone  with angle $\theta \leq \cos^{-1}(1/\sqrt d)$.
This is in stark contrast with the cases of the Laplacian and the fractional Laplacian, for which the  Varopoulos-type Dirichlet heat kernel estimates hold in any Lipschitz open set (see \cite{V} and \cite{BGR1}).  Since a $C^1$ open set is a Lipschitz open set satisfying the interior cone condition with common angle $\theta\in (\cos^{-1}(1/\sqrt d), \pi)$  (see \cite{CW2}), as a corollary of Theorem \ref{T3}, the Varopoulos-type Dirichlet heat kernel estimates holds for  $\Delta+\Delta^{\alpha/2}$ in any $C^1$ open set.

(ii) The scale invariant  boundary Harnack principle is the main tool   in the study of the Dirichlet heat kernel estimates for  pure jump L\'evy processes  in $\kappa$-fat  open sets  (see \cite{BGR1, CKS6}).
Due to the fact mentioned above that the scale invariant boundary Harnack principle for $X$ may fail to hold in some Lipschitz open sets,  the method of the scale invariant boundary Harnack principle does not work in the proof of Theorem \ref{T2}.
In \cite{ChoKSV} for  purely non-local operators with killing potentials in $\kappa$-fat open sets, the Varopoulos-type Dirichlet heat kernel estimates for the operators are obtained by Dynkin's formula and some condition on the  operator  instead of the scale invariant boundary Harnack principle. While the proof in \cite{ChoKSV} relies crucially on the pure jump property of the corresponding jump processes, its argument  could not be applied  to our setting.
New ideas  are required to obtain the two-sided Dirichlet heat kernel estimates for  L\'evy processes with Gaussian components  in  Lipschitz open sets.

 (iii) For the purely non-local operators  in the literature \cite{BGR1, CKS6, ChoKSV},  the survival probability  for the pure jump processes in $\kappa$-fat open sets  is shown to be  comparable to the exit time for the  process  multiplied by $t^{-1}$.
Distinct from  purely non-local operators, the key observation in this paper is that the survival probability $\P_x(\tau_D>t)$ for the process $X$ in a Lipschitz open set $D$ is  comparable to  the Green function for $X$ in a small open set near the boundary multiplied by $t^{(d-2)/2}$ (see Proposition \ref{L:4.5}).
This property  is also the key step in the proof of Theorem \ref{T3}.

(iv) The method in  Varopoulos \cite{V} for  differential operators in Lipschitz open sets is mainly analytic. The upper bound estimates of Dirichlet heat kernel  in \cite{V}  is established mainly by the Gaussian mass escape  and the iteration method in  \cite{U}.  The lower bound estimates of Dirichlet heat kernel  in \cite{V}  is obtained by the scale invariant parabolic boundary Harnack principle.  In our setting, the mass escape  for the process $X$ in large scale is dominated  by $\frac{t}{r^d \phi (r)}$  (see Theorem \ref{T1}). This differs from the Gaussian case in \cite{V}.
Moreover, Theorems \ref{T2} and \ref{T3} show that the boundary term of $p_D(t, x, y)$ in large scale $|x-y|$ may not be comparable to the survival probability. This reflects the  difference between the operator $\L$ in this paper and the differential operator.

\end{remark}

In this paper, different from \cite{V} for differential operators, the method  is mainly probabilistic. Instead of the scale invariant boundary Harnack principle, the main tools  in this paper are the Carleson estimate for harmonic function for $X$ in Lipschitz open sets, the parabolic Harnack inequality and one type of parabolic Carleson estimate with respect to the time-space process of $X$ in Lipschitz open sets.
Some probabilistic argument  are based on \cite{BGR1, CKS6}.
The doubling property for the survival probability  $\P_x(\tau_D>t)$ is inspired by the method of the parabolic Harnack inequality and the parabolic Carleson estimate in \cite{V}.
The strategy of the proof is as follows.

For the upper bound estimates of $p_D(t, x, y)$ in Theorem \ref{T2}, by the strong Markov property of $X$ and the probabilistic formula of $p_D(t, x, y)$ (see Lemma \ref{L:4.1} below),  the key ingredient  is to show the exit distribution for $X$ from a small open subset of $D$ near the boundary is bounded from above by the survival probability $\P_x(\tau_D>t)$.
Different from  the method of the scale invariant boundary Harnack principle in  \cite{BGR1, CKS6} and \cite{ChoKSV} for purely non-local operators, we mainly use Dynkin's formula and applying the Carleson estimate of harmonic function with respect to $X$ in Lipschitz open sets  to derive the upper bound estimates of the  exit distribution  for $X$ from a small open subset of $D$ (see Proposition \ref{P:2.1}). By this estimate combining with the hitting probability in Lemma \ref{L:2.6} and the strong Markov property of $X,$  it can be obtained that  the exit distribution for $X$ from a small open subset of $D$ near the boundary is bounded from above by the survival probability (see Proposition \ref{P:2.2}).

For the lower bound estimates of Dirichlet heat kernel,
 the  scale invariant parabolic boundary Harnack inequality is the main tool in \cite{V} for the  differential operators in  Lipschitz open sets, especially in bounded Lipschitz open sets. The L\'evy system formula for jump processes and the scale invariant boundary Harnack principle are mainly used in \cite{BGR1, CKS6} for purely non-local operators.
 In our setting, to obtain  the lower bound estimates of $p_D(t, x, y)$ in Theorem \ref{T2},  we  use instead the parabolic Carleson estimate for a class of nonnegative parabolic functions with respect to the time-space process $Z_t=(V_0-t, X_t)$ in cylindrical set $(0, 4r^2)\times U$  that vanish on $(0, 4r^2)\times U^c$ for $r\in (0, T)$ and a small Lipschitz subset $U$ of $D$ near the boundary(see Proposition \ref{P:3.2}) and the parabolic Harnack inequality for $Z_t=(V_0-t, X_t)$ established in \cite{CK3}.
On the other hand, the  parabolic Harnack inequality and the parabolic Carleson estimate for differential operators are mainly used in the proof of the doubling property for the survival probability of diffusion in \cite[Section 2]{V}.
 The  doubling property for the survival probability of the process $X$ (see Proposition \ref{P:3.4}) is  inspired by the idea in \cite{V}.
While our proof is mainly probabilistic, the argument  differs from \cite{V}.

 The necessary and sufficient conditions for the Varopoulos-type Dirichlet heat kernel estimates  for the process $X$ in a Lipschitz open set $D$ are given in Theorem \ref{T3}.  As mentioned above, the comparability between the surviving probability $\P_x(\tau_D>t)$ and the Green function in a small open set near the boundary multiplied by $t^{(d-2)/2}$  in Proposition \ref{L:4.5} is the crucial step in the proof of Theorem \ref{T3}.
This comparability is mainly obtained by the upper bound estimates of Dirichlet heat kernel in Theorem \ref{T2}, the exit distribution estimates in Proposition \ref{P:2.1} and the   hitting probability in Lemma \ref{L:2.6}.
The equivalence of \eqref{e:1.10} and the scale invariant boundary Harnack principle \eqref{e:1.13} in Theorem \ref{T3} is obtained by the result of the scale invariant boundary Harnack principle for discontinuous L\'evy processes with Gaussian components in \cite{CW2}.

\medskip

In the last part, as an application of Theorem \ref{T3},   the explicit two-sided Dirichlet heat kernel estimates for $X$ in $C^{1, 1}$ open sets  are obtained.
This result recovers the case  in \cite{CKS7} for subordinate Brownian motion with Gaussian components   in $C^{1, 1}$ open sets in $\R^d (d\geq 3)$  when the scaling order of the pure jump part of the subordinator is strictly between $0$ and $1$.

For an open set $D$ and a positive constant $c,$ we define
$$h_{D, c}(t, x, y):=\left(1\wedge \dfrac{\delta_D(x)}{\sqrt{t}}\right)\left(1\wedge \dfrac{\delta_D(y)}{\sqrt{t}}\right) \left(t^{-d/2}\wedge \left( t^{-d/2}e^{-c|x-y|^2/t}+ t^{-d/2} \wedge
tj(|x-y|)\right) \right).$$

\begin{cor}\label{C2}
Suppose that $D$ is a $C^{1, 1}$ open set in $\R^d$ for $d\geq 3$ with characteristics $(R_0, \Lambda_0)$ and the path distance in each connected component of $D$ is comparable to the Euclidean distance with characteristic $\chi_1.$
For each $T>0,$ there exist positive constants $C_1=C_1(d, \phi, R_0, \Lambda_0, \chi_1, T)$ and $C_k=C_k(d, \phi, \chi_1), k=2, 3$ such that for  all $x, y\in D$ and $t\in (0, T),$
\begin{equation}\label{e:1.14}
C_1^{-1}h_{D, C_2}(t, x, y)
\leq p_D(t, x, y)\leq C_1h_{D,C_3}(t, x, y).
\end{equation}
Suppose in addition that $D$ is bounded. Then there exists a constant $C_4=C_4(d, \phi, R_0, \Lambda_0, {\rm diam}(D))$ such that for $t\geq 3,$
\begin{equation}\label{e:1.15}
C_4^{-1}e^{-\lambda_1 t}\delta_D(x)\delta_D(y)\leq p_D(t, x, y)\leq C_4e^{-\lambda_1 t}\delta_D(x)\delta_D(y),
\end{equation}
where $-\lambda_1<0$ is the largest eigenvalue of the generator of $X^D.$
\end{cor}

The organization of this paper is as follows.  Some preliminary results are given in Section 2. In Section 3, we first obtain the upper bound estimates of the exit distribution for the process $X$ from a small open subset  near the boundary of a Lipschitz open set $D$ in Proposition \ref{P:2.1}.  By this estimate combining with the hitting probability in Lemma \ref{L:2.6} and the strong Markov property of $X,$  the exit distribution for $X$ from a small open subset of $D$ near the boundary is shown to be bounded from above by the survival probability in Proposition \ref{P:2.2}.  In Section 4, we first establish one type of parabolic Carleson estimate for nonnegative parabolic functions with respect to  $Z_t=(V_0-t, X_t)$ in cylindrical set $(0, 4r^2)\times U$  that vanish on $(0, 4r^2)\times U^c$ for $r\in (0, T)$ and a small Lipschitz subset $U$ of $D$ near the boundary.  By this parabolic Carleson estimate and the parabolic Harnack inequality established in \cite{CK3},  the  doubling property of the survival probability for $X$ is obtained in Proposition \ref{P:3.4}. Section 5 is devoted to prove  Theorem \ref{T2}.  In Section 6, we prove Theorem \ref{T3} and Corollary \ref{C2}.

\section {Preliminaries}

In this section,  some preliminary results will be given.
For each open set $B,$ denote by $\tau_B$ the first exit time of $X$ from $B.$
By virtue of Lemmas 2.3 and 2.4 in \cite{CK3}, the following result holds.

\begin{lem}\label{L:2.2}
Let $r_0>0.$ There exists   $C=C(d, \phi, r_0)>1$ such that for any $x_0\in\R^d$ and  $r\in (0, r_0),$
$$
C^{-1}r^2\leq \E_{x_0}\tau_{B(x_0, r)}  \leq Cr^2.
$$
\end{lem}

Let
  $Z=\{Z_t:=(V_0-t, X_t), t\geq 0\}$
 denote the
time-space process of $X$.
Denote by $\P_{(t, x)}$ the law of the time-space process $Z$ starting from $(t, x).$
Let $D$ be an open set.
We say that  a non-negative real valued Borel  measurable function
$h(t,x)$ on $[0, \infty)\times \R^d$ is {\it parabolic}
(or {\it caloric})
on $Q=(a,b)\times D\subset [0,\infty)\times\R^d$
if
for every relatively compact open subset $Q_1$ of $Q$,
$$
h(t, x)=\E_{(t,x)} [h (Z_{\tau_{Q_1}})]
$$
for every $(t, x)\in Q_1$,
where
$\tau_{Q_1}=\inf\{s> 0: \, Z_s\notin Q_1\}$.

The following parabolic Harnack inequality  for the time-space process $Z=\{Z_t:=(V_0-t, X_t), t\geq 0\}$ is established in Theorem 1.3 of \cite{CK3}.

\begin{prp}[Parabolic Harnack inequality]\label{P:3.1}
For every $\delta\in (0, 1)$, there
exist positive constants $c_k=c_k(d, \phi, \delta), k=1, 2$
such that for every $x\in \R^d$, $t_0\ge 0$,
$0<R\leq c_1$
and every non-negative function $u$ on $[0, \infty)\times \R^d$
that is parabolic with respect to $Z$
on $(t_0,t_0+6\delta R^2 )\times B(x,4R)$,
$$
\sup_{(t_1,y_1)\in Q_-}u(t_1,y_1)
\le c_2 \,
\inf_{(t_2,y_2)\in Q_+}u(t_2,y_2),
$$
where $Q_-=(t_0+\delta R^2,t_0+2\delta R^2)\times B(x, R)$ and
$Q_+=(t_0+3\delta R^2,t_0+ 4\delta R^2)\times B(x, R)$.
\end{prp}

A real-valued function $h$ defined on $\R^d$ is said to be
 {\it harmonic} in an open set  $U$ with respect to $X$
if for every open set $B$  whose closure is a compact subset of $U,$
   $\E_x |h(X_{\tau_B})|<\infty$ and
 $h(x)=\E_x h(X_{\tau_B})$ for  each  $x\in B.$
  In particular, we say $h$  is   regular harmonic in   $U$ with respect to $X$
if    $\E_x |h(X_{\tau_U})|<\infty$ and
$h(x)=\E_x h(X_{\tau_U})$  for each  $x\in U.$

As a particular case of Proposition \ref{P:3.1}, the elliptic Harnack inequality holds for $X.$ That is, there
exists a constant $c=c(d, \phi)>0$
such that for every $x\in \R^d, 0<r\leq 1$
and every non-negative function $h$ on $\R^d$
that is  harmonic with respect to $X$
on $B(x,4r)$,
\begin{equation}\label{e:2.5}
h(y_1)\leq ch(y_2), \quad y_1, y_2\in B(x, r).
\end{equation}

\smallskip

It is known that any discontinuous Hunt process admits a L\'evy system that describes how the process jumps.
By  \cite[Lemma 4.2]{CK3}, for any  nonnegative function $f$ on $\R_+\times
\R^d\times \R^d$ vanishing along the diagonal of $\R^d \times \R^d$, for any stopping time $T$
with respect to the minimal admissible  augmented  filtration generated by $X$ and $x\in \R^d$,
\begin{equation}\label{e:2.2}
\E_x \left[ \sum_{s\leq T}f(s, X_{s-}, X_s); X_{s-}\neq X_s \right]
=\E_x \left[ \int_0^T\int_{\R^d}
f(s,X_s,y) j(|X_s-y|)\,dy\,ds \right].
\end{equation}

\begin{lem}\label{L:2.7}
For each $b>0,$ there exists $c=c(d, \phi, b)\in (0, 1)$ such that for any $x_0\in \R^d$ and $r\in (0, 1),$
$$\P_{x_0}(\tau_{B(x_0, r)}\geq b r^2)\geq c.$$
\end{lem}

\pf By  \cite[Theorem 3.1]{CK3}, there exists $c_0=c_0(d, \phi)>0$ such that for any $x_0\in\R^d$ and $r\in (0, 1),$
\begin{equation}\label{e:2.9}
p_{B(x_0, r)}(r^2, x, y)\geq c_0r^{-d}, \quad x, y\in B(x_0, r/2).
\end{equation}
If $b\in (0, 1],$ then by \eqref{e:2.9},
$$\P_{x_0}(\tau_{B(x_0, r)}\geq b r^2)\geq \P_{x_0}(\tau_{B(x_0, r)}\geq r^2)\geq \int_{B(x_0, r/2)} p_{B(x_0, r)}(r^2, x_0, y)\,dy \geq c_1.$$
If $b>1,$ let $[b]$ be the integer part of $b$ and  $m=[b]+1,$  then by  \eqref{e:2.9},
$$\begin{aligned}
&\P_{x_0}(\tau_{B(x_0, r)}\geq b r^2)\geq \P_{x_0}(\tau_{B(x_0, r)}\geq m r^2)\geq \int_{B(x_0, r/2)} p_{B(x_0, r)}(mr^2, x_0, y)\,dy\\
&\geq  \int_{B(y_{m-1}, \frac{r}{4m})}\cdots \int_{B(y_1, \frac{r}{4m})}\int_{B(x_0, \frac{r}{4m})} p_{B(x_0, r)}(r^2, x_0, y_1)p_{B(x_0, r)}(r^2, y_1, y_2)\\
&\qquad\cdots p_{B(x_0, r)}(r^2, y_{m-1}, y_m)\,dy_1\cdots dy_m\\
&\geq  \int_{B(y_{m-1}, \frac{r}{4m})}\cdots \int_{B(y_1, \frac{r}{4m})}\int_{B(x_0, \frac{r}{4m})} (c_0r^{-d})^{m}\,dy_1\cdots dy_m\\
&\geq  c_2,
\end{aligned}$$
where  $c_2=c_2(d, \phi, m),$ the second line from the bottom is due to  $|y_i-x_0|\leq \sum_{k=2}^i |y_k-y_{k-1}|+|y_1-x_0|\leq \frac{ir}{4m}\leq r/2$ for $i=1, \cdots, m$  and \eqref{e:2.9}. The proof is complete.

\qed

\smallskip

Recall that for each open set $B,$ we denote by $G_B(x, y)$ the Green function of $X$ killed upon exiting $B.$

\begin{lem}\label{L:2.4'}
There exist $C_k=C_k(d, \phi)>0, k=1, 2$ such that for any $x\in \R^d$ and $r\in (0, 1),$
$$G_{B(x, r)}(x, y)\geq C_1r^{2-d} \quad \mbox{for } \quad y\in B(x, r/4)$$
and
$$G_{B(x, r)}(x, y)\leq C_2r^{2-d} \quad \mbox{for }\quad   y\in B(x, r/2) \setminus B(x, r/4).$$
\end{lem}

\pf By \eqref{e:2.9}, there exists $c_1=c_1(d, \phi)>0$ such that for $x\in\R^d$ and $y\in B(x, r/4),$
$$G_{B(x, r)}(x, y)=\int_0^\infty p_{B(x, r)}(t, x, y)\,dt\geq \int_{r^2/2}^{r^2} p_{B(x, r)}(t, x, y)\,dt\geq c_1r^{2-d}.$$
Note that $G_{B(x, r)}(x, \cdot)$ is harmonic in $B(x, r)\setminus \{x\}.$ By the elliptic Harnack inequality \eqref{e:2.5} and the standard chain argument, there exists $c_2=c_2(d, \phi)>0$ such that
$$\sup_{y\in B(x, r/2) \setminus B(x, r/4)}G_{B(x, r)}(x, y)\leq c_2\inf_{u\in B(x, r/2) \setminus B(x, r/4)}G_{B(x, r)}(x, u).$$
Hence by Lemma \ref{L:2.2}, there exists $c_3=c_3(d, \phi)>0$ such that for $y\in B(x, r/2) \setminus B(x, r/4),$
$$G_{B(x, r)}(x, y)\leq c_2\dfrac{1}{r^d} \int_{B(x, r/2) \setminus B(x, r/4)} G_{B(x, r)}(x, u)\,du\leq c_2\dfrac{1}{r^d}\E_x \tau_{B(x, r)}\leq c_3r^{2-d}.$$
The proof is complete.

\qed

\section{ Exit distribution  and survival probability estimates}

In this section, we first obtain the upper bound estimates of the exit distribution for  $X$ from a small open subset  near the boundary of a Lipschitz open set $D.$ By this estimate, we show  the exit distribution above  is  bounded from above by the survival probability. Furthermore  the comparability between the survival probability of $X$ from $D$ and that from a small open subset of $D$ near the boundary is obtained.

\begin{defn}\label{D1}

An open set $D$ in $\R^d$ (when $d\geq 2$) is said to be Lipschitz if there exist a localization radius
$R_0 >0$ and a constant $\Lambda_0>0$ such that for every $z\in\partial D,$ there exist a Lipschitz function
$\Gamma=\Gamma_z: \R^{d-1}\rightarrow \R$ satisfying
$$
\Gamma(0)=0, \quad
|\Gamma(x)-\Gamma(y)|\leq \Lambda_0 |x-y|
$$
 and an orthonormal coordinate system
$CS_z: y=(y^1, \cdots, y^{d-1}, y^d)=:(\wt {y}, y^d)\in \R^{d-1}\times \R$ with its origin at $z$ such that
$$
B(z, R_0)\cap D=\{y=(\wt {y}, y^d)\in B(0,R_0) \hbox{ in } CS_z: y^d>\Gamma(\wt {y})\}.
$$
The pair $(R_0, \Lambda_0)$ is called the characteristics of the Lipschitz open set $D.$
\end{defn}

Without loss of generality, throughout this paper we assume $R_0\leq 1$ and $\Lambda_0\geq 1.$
For any $x\in D,$ let $\delta_D(x)$ denote the Euclidean distance between $x$ and $D^c$.
Let $D$ be a Lipschitz open set with characteristics $(R_0, \Lambda_0).$
Recall  that
there exists   $ \kappa=\kappa(d, \Lambda_0)\in (0, 1/4)$ such that for $r\in (0, R_0)$
and $z\in\partial D,$
\begin{equation}\label{e:1.7}
\hbox{there exists } z_r\in D\cap \partial B(z, r)  \hbox{ with }
\kappa r\leq \delta_D(z_r)< r.
 \end{equation}
 In the remainder of this paper, we always use $\kappa$ to denote the constant in \eqref{e:1.7}.

\begin{defn}\label{D2}
We say that the path distance in a domain (i.e. connected open set) $U$ is comparable to the Euclidean distance with characteristic $\chi_1$
if for every $x$ and $y$ in $U,$ there is a rectifiable curve $l$ in $U$ which connects $x$ to $y$ such that the length of $l$ is less than or equal to $\chi_1|x-y|.$ Clearly such a property holds for all bounded Lipschitz domains, Lipschitz domains with compact complements and domain consisting of all the points above the graph of  Lipschitz function.
\end{defn}

 \medskip

 For   $\theta\in (0, \pi),$ define
\begin{equation}\label{e:cone}
C_\theta:=\{x=(x^1, \ldots, x^d) \in \R^d: x^d>|x|\cos\theta\},
\end{equation}
 which is a    circular cone with apex $2\theta$ and axis along the positive  $x^d$-direction.
 We    call $\theta$ the angle of the cone $C_\theta.$

\begin{lem}\label{L:2.1}
There exists a positive constant $\delta=\delta(d, \phi)\in (0, 1)$ such that for any $ r\in (0, \delta]$ and any circular cone $C_\theta$ with the angle $\theta,$
\begin{equation}\label{e:2.1'}
\P_0(X_{\tau_{B(0, r)}}\in C_\theta\cap B(0, 2r))\geq \dfrac{\theta}{2\pi}.
\end{equation}
\end{lem}

\pf  We first prove that there exists a constant $\delta\in (0, 1)$ such that for any $r\in (0, \delta),$
\begin{equation}\label{e:2.1}
\P_0(X_{\tau_{B(0, r)}}\in B(0, 2r))\geq \dfrac{1}{2}.
\end{equation}
In fact, by the L\'evy system formula \eqref{e:2.2}  and Lemma \ref{L:2.2}, for $r\in (0, 1),$
$$\begin{aligned}
\P_0(X_{\tau_{B(0, r)}}\in B^c(0, 2r))&=\E_0 \int_0^{\tau_{B(0, r)}}\int_{B^c(0, 2r)} \dfrac{1}{|X_s-y|^d \phi(|X_s-y|)} \,dy\,ds\\
&\leq c_1\E_0 \int_0^{\tau_{B(0, r)}}\int_{B^c(0, 2r)} \dfrac{1}{|y|^d \phi(|y|)} \,dy\,ds\\
&\leq c_2\E_0 \int_0^{\tau_{B(0, r)}}\int_{B^c(0, 2r)} \left[\dfrac{1}{|y|^{d+\beta_2}}+ \dfrac{1}{|y|^{d+\beta_1}}\right] \,dy\,ds\\
&\leq c_3\E_0 \tau_{B(0, r)}r^{-\beta_2}\\
&\leq c_4r^{2-\beta_2},
\end{aligned}$$
where $c_k=c_k(d, \phi), k=1, \cdots, 4$ and  we used  \eqref{e:1.1} and \eqref{e:1.3}-\eqref{e:1.4} in the second and third lines.
Hence we can take  $\delta\in (0, 1)$ small enough such that for any $r\in (0, \delta),$
$$\P_0(X_{\tau_{B(0, r)}}\in B^c(0, 2r))\leq \dfrac{1}{2}.$$
That is, \eqref{e:2.1} holds.
Let $T$ denote any rotation operator in $\R^d$ and set the process $\tilde X=TX.$
Since the distribution of $X$ is  rotationally symmetric, the process $\tilde X$ and $X$ are identical in the sense of the distribution.
 We have for $r\in (0, \delta),$
$$\begin{aligned}
\P_0(X_{\tau_{B(0, r)}}\in C_\theta\cap B(0, 2r))&=\P_0(\tilde X_{\tau_{B(0, r)}}\in T\circ (C_\theta\cap B(0, 2r)))\\
&=\P_0(X_{\tau_{B(0, r)}}\in T\circ (C_\theta\cap B(0, 2r)))\\
&=\P_0(X_{\tau_{B(0, r)}}\in (T\circ C_\theta)\cap B(0, 2r)).
\end{aligned}$$
Note that for any rotation $T$ in $\R^d,$  $T\circ C_\theta$ is also a cone with the angle $\theta$  but  with the axis along the rotation $T$ of the positive  $x^d$-direction.
Hence
\begin{equation}\label{e:2.5'}
\P_0(X_{\tau_{B(0, r)}}\in C_\theta\cap B(0, 2r))=\dfrac{\theta}{\pi}\cdot \P_0(X_{\tau_{B(0, r)}}\in B(0, 2r)).
\end{equation}
Thus the desired conclusion \eqref{e:2.1'} can be obtained by  \eqref{e:2.1} and \eqref{e:2.5'}.

\qed

 By an argument similar to that of \cite[Lemma 4.2]{CKSV}, we have the following result.

\begin{lem}\label{L:2.1'}
Let $D$ be an open set in $\R^d$ and $U$ and $V$ be bounded open sets
with $V \subset \overline{V} \subset U$ and $ D \cap V \not= \emptyset$.
Suppose $u$ is a nonnegative function in $\R^d$ that is
harmonic in $D\cap U$ with respect to $X$ and vanishes
continuously on $D^c\cap U$. Then $u$ is regular harmonic in
$D\cap V$ with respect to $X$, i.e.,
\begin{equation}
u(x)=\E_x\left[ u(X_{\tau_{D\cap V}})\right] \qquad \hbox{
for all }x\in D\cap V\, .
\end{equation}
\end{lem}

\begin{prp}[Carleson estimate]\label{L:2.4}
Suppose $D$ is a Lipschitz open set with characteristics $(R_0, \Lambda_0)$.
There exists a constant  $C=C(d, \phi, \Lambda_0)>0$   such that  for any $z_0\in\partial D, r\in (0, R_0/4)$ and any nonnegative harmonic function $h$
on $D\cap B(z_0, 2r)$  vanishing continuously on  $D^c\cap B(z_0, 2r),$
\begin{equation}\label{e:2.9nn}
h(x)\leq Ch(x_0)  \quad \hbox{for } x\in D\cap B(z_0, r),
\end{equation}
where $x_0$ is a point of $D\cap \partial B(z_0, r)$ with $\kappa r\leq\delta_D(x_0)< r.$
\end{prp}

 \pf   Let $\delta$ be the constant in Lemma \ref{L:2.1}. Let $z_0\in \partial D$ and  $r\in (0, R_0/4).$  By the  elliptic Harnack inequality \eqref{e:2.5} and the standard chain argument, it suffices to prove \eqref{e:2.9nn} for
 $x\in D\cap B(z_0, r\wedge \frac{\delta}{2}).$
  Let   $x\in D\cap B(z_0, r\wedge \frac{\delta}{2})$ and
  $z_x$ be a point of $\partial D$ such that $|x-z_x|=\delta_D(x).$
 Since the Lipschitz open set $D$ satisfies the exterior cone condition, there exist an angle $\theta=\theta(d, \Lambda_0)\in (0, \pi)$ and a cone $C_\theta(z_x)$ that is conjugate to $C_\theta$ with  angle $\theta$  and vertex at $z_x$ such that
 $C_\theta(z_x)\cap B(z_x, R_0)\subset D^c.$
 Note that $D\cap B(z_x, 2\delta_D(x))\subset D\cap B(x, 4\delta_D(x)).$
Hence
 \begin{equation}\label{e:2.11'}
 \begin{aligned}
\P_x(X_{\tau_{D\cap B(x, 4\delta_D(x))}}\in D)&\leq \P_x(X_{\tau_{D\cap B(z_x, 2\delta_D(x))}}\in D)\\
&\leq \P_x(X_{\tau_{B(z_x, 2\delta_D(x))}}\in D)\\
&=1- \P_x(X_{\tau_{B(z_x, 2\delta_D(x))}}\in D^c)\\
&\leq 1-\P_x(X_{\tau_{B(z_x, 2\delta_D(x))}}\in C_\theta(z_x)\cap B(z_x, 4\delta_D(x))).
 \end{aligned}\end{equation}
 Note that $h(y):=\P_y(X_{\tau_{B(z_x, 2\delta_D(x))}}\in C_\theta(z_x)\cap B(z_x, 4\delta_D(x)))$ is a harmonic function with respect to $X$ in $B(z_x, 2\delta_D(x)).$ By the uniform Harnack inequality \eqref{e:2.5}, there exists $c_1=c_1(d, \phi)\in (0, 1)$ such that
\begin{equation}\label{e:2.9''}
\P_x(X_{\tau_{B(z_x, 2\delta_D(x))}}\in C_\theta(z_x)\cap B(z_x, 4\delta_D(x)))
\geq c_1\P_{z_x}(X_{\tau_{B(z_x, 2\delta_D(x))}}\in C_\theta(z_x)\cap B(z_x, 4\delta_D(x))).
\end{equation}
 By Lemma \ref{L:2.1},  we have
\begin{equation}\label{e:2.8n}
\P_{z_x}(X_{\tau_{B(z_x, 2\delta_D(x))}}\in C_\theta(z_x)\cap B(z_x, 4\delta_D(x)))\geq \dfrac{\theta}{2\pi}.
\end{equation}
Hence by combining \eqref{e:2.9''} and \eqref{e:2.8n},
\begin{equation}\label{e:2.12}
\P_x(X_{\tau_{B(z_x, 2\delta_D(x))}}\in C_\theta(z_x)\cap B(z_x, 4\delta_D(x)))
\geq c_1\dfrac{\theta}{2\pi}.
\end{equation}
Let $c_2:=c_1\dfrac{\theta}{2\pi}\in (0, 1).$
Thus by \eqref{e:2.11'} and \eqref{e:2.12},
\begin{equation}\label{e:2.3}
\P_x(X_{\tau_{D\cap B(x, 4\delta_D(x))}}\in D)\leq 1-c_2\in (0, 1).
\end{equation}
By the condition \eqref{e:1.1}, it is easy to see that $j(|z|)\asymp j(2|z|)$ for $|z|\leq 1$ and $j(|z|)\asymp j(|z|+1)$ for $|z|>1.$
Note that the constant $c_2$ in \eqref{e:2.3} only depends on $(d, \phi, \Lambda_0).$ Moreover   the elliptic Harnack inequality \eqref{e:2.5} and Lemma \ref{L:2.2} hold with the comparison constants  only depending on $(d, \phi)$.
By  an argument similar to  that of \cite[Theorem 4.3]{CKSV}  with \eqref{e:2.3} in place of \cite[Lemma 4.1]{CKSV},
 the desired conclusion \eqref{e:2.9nn} can be obtained.
\qed

\medskip

Let $D$ be a Lipschitz open set  with characteristics $(R_0, \Lambda_0).$
It is known that (see e.g.  \cite[P4 of Section 5]{Bo}) there exists $L=L(\Lambda_0, d)>1$
such that for any $z\in\partial D$ and $r\in (0, R_0),$  there is a Lipschitz connected open set  $U_r(z)\subset D$   with characteristics $(r/L, L\Lambda_0)$ such that $D\cap B(z, r/L)\subset U_r(z)\subset D\cap B(z, r)$.
It follows from Remark 5 in \cite{Bo} ( see also in \cite[page 136]{ChZ}) that we can choose $U_r(z)$ a cylinder set in the coordinate system $CS_z$ as the intersection of $D$ and a cylinder along the $y^d$ axis.
 In detail, we choose
$U_r(z)=\{y=(\tilde y, y^d)\: \mbox{in} \: CS_z: |\tilde y|<ar, \: \Gamma_z(\tilde y)<y^d<\Gamma_z(\tilde y)+br\}$ for some positive constants $a$ and $b$  such that $D\cap B(z, r/L)\subset U_r(z)\subset D\cap B(z, r)$,
where $\Gamma_z$ is the Lipschitz function with Lipschitz constant $\Lambda_0$ defined in Definition \ref{D1}.
In the following, we always use $L$ and $U_r(z)$ to denote such constant and such a local Lipschitz open set.

For each $z\in \partial D$ and $r\in (0, R_0),$ as noted above,  $U_r(z)$  is above the graph of the Lipschitz function $\Gamma_z$ in the coordinate system $CS_z.$
By an argument similar to that of   case (b) in  the proof of \cite[Lemma 6.6(iii)]{ChZ}, for any given $M>0$ and $x, y\in U_r(z)$ such that $|x-y|\geq M(\delta_{U_r(z)}(x)\wedge \delta_{U_r(z)}(y)),$ there exist a finite number  of balls $B(a_i, r_i)(0\leq i\leq N)$ in $U_r(z)$ such that $a_0=x, a_N=y$ and $B(a_i, \frac{r_i}{2})\cap B(a_{i+1}, \frac{r_{i+1}}{2})\neq \varnothing$ for $0\leq i\leq N-1.$ Here the number $N$ depends only on $ M$ and the Lipschitz constant $\Lambda_0.$

\begin{prp}\label{P:2.1}
Suppose $D$ is a Lipschitz open set with characteristics $(R_0, \Lambda_0).$
There exists $C=C(d, \phi, \Lambda_0)>0$ such that for any $z_0\in \partial D$ and $r\in (0, R_0),$
$$\P_x(X_{\tau_{U_r(z_0)}}\in D)\leq C\dfrac{G_{U_r(z_0)}(x, y_0)}{r^{2-d}} \qquad \hbox{for} \quad x\in D\cap B(z_0, r/(8L)),$$
where  $y_0$ is a point of  $D\cap \partial B(z_0, r/(2L))$ with $\delta_D(y_0)\geq \kappa r/(2L).$
\end{prp}

\pf Let $\varphi$ be a non-negative $C^\infty(\R^d)$ function such that $\varphi=0$ on $B(0, 3/(4L))$ and $\varphi=1$ on $\R^d\setminus B(0, 1/L).$
Let $z_0\in \partial D$ and $r\in (0, R_0).$
For each $n \geq 1,$ let $\varphi^{(n)}_{z_0, r}$ be a $C_c^\infty(\R^d)$ function taking values in $[0, 1]$ such that $\varphi^{(n)}_{z_0, r}(y)=\varphi((y-z_0)/r)$ for $y\in B(z_0, n).$
By \eqref{e:L} and the fact that $j(z)=j(|z|)$ is rotationally symmetric and satisfies $\int_{\R^d} (1\wedge |z|^2) j(|z|)<\infty$, for any $r>0$ and $f\in C_c^\infty(\R^d),$
\begin{equation}\label{e:2.17n}\begin{aligned}
\L f(y)&=\Delta f(y)+\lim_{\ee\rightarrow 0}\int_{\{z\in \R^d: |z|>\ee\}} (f(y+z)-f(y))j(z)\,dz\\
 &=\Delta f(y)+\int_{|z|\leq  r } \left(f(y+z)-f(y)-\nabla f(y) \cdot z\right) j(z)\,dz\\
 &\qquad +\int_{|z|> r } \left(f(y+z)-f(y)\right) j(z)\,dz, \quad y\in\R^d.
\end{aligned}\end{equation}
Note that for each $n\geq 2$ and $r\in (0, R_0)\subset (0, 1),$ $\varphi^{(n)}_{z_0, r}(y)=\varphi((y-z_0)/r)$ for $y\in B(z_0, 2r)\subset B(z_0, 2).$
Hence for each $n\geq 2$ and $r\in (0, R_0),$
$\dfrac{\partial^2}{\partial y_iy_j} \varphi^{(n)}_{z_0, r}(y)=r^{-2}\dfrac{\partial^2 \varphi}{\partial y_iy_j}  ((y-z_0)/r) $  and $\nabla   \varphi^{(n)}_{z_0, r} (y) = r^{-1} \nabla \varphi ((y-z_0)/r)$ for $y\in B(z_0, 2r).$
By \eqref{e:2.17n}, \eqref{e:1.3} and \eqref{e:1.4}, there exist $c_k=c_k(d, \phi), k=1, 2$  such that for any $n\geq 2, r\in (0, R_0)$ and $y\in B(z_0, r),$
\begin{equation}\label{e:2.7}\begin{aligned}
&|\L \varphi^{(n)}_{z_0,r}(y)|\\
\leq &\left|\Delta \varphi^{(n)}_{z_0, r}(y)\right|+\int_{|z|\leq  r } \left|\varphi^{(n)}_{z_0, r}(y+z)-\varphi^{(n)}_{z_0,r}(y)-\nabla \varphi^{(n)}_{z_0, r} (y) \cdot z\right| j(|z|)\,dz\\
&+\int_{|z|> r } |\varphi^{(n)}_{z_0, r}(y+z)-\varphi^{(n)}_{z_0, r}(y)| j(|z|)\,dz\\
\leq &r^{-2}\left( \| D^2 \varphi\|_\infty +c_1\| D^2 \varphi\|_\infty \ \int_{|z|\leq  r } |z|^2 |z|^{-(d+\beta_2)}\,dz\right)\\
&+ c_1\int_{r<|z|\leq 1 }  |z|^{-(d+\beta_2)}\,dz +c_1\int_{|z|> 1 }  |z|^{-(d+\beta_1)}\,dz\\
\leq & c_2 r^{-2}.
\end{aligned}\end{equation}

Let $r\in (0, R_0)$ and fix $x\in D\cap B(z_0, r/(8L)).$
Let $y_0$ be a point of  $D\cap \partial B(z_0, r/(2L))$ with $\delta_D(y_0)\geq \kappa r/(2L).$
For the simplicity of notation,  we denote $D_r(z_0):=D\cap B(z_0, r).$
Note that $D_{r/L}(z_0) \subset U_r(z_0) \subset D_r(z_0),$  thus $\delta_{U_r(z_0)}(y_0)\geq \kappa r/(2L).$
 For each $ y\in U_r(z_0),$ let $z_y\in \partial U_r(z_0)$ such that
$|y-z_y|=\delta_{U_r(z_0)}(y).$ Note that $G_{U_r(z_0)}(x, \cdot)$ is harmonic in $U_r(z_0)\setminus \{x\}$ and vanishes on $U_r(z_0)^c.$
Then for each $ y\in U_r(z_0) \setminus D_{r/(2L)}(z_0)$ with $\delta_{U_r(z_0)}(y)\leq r/(16L),$
$G_{U_r(z_0)}(x, \cdot)$ is harmonic in $U_r(z_0)\cap B(z_y, r/(8L))$ and vanishes on $U_r(z_0)^c.$
By applying the Carleson estimate of Proposition \ref{L:2.4} on $G_{U_r(z_0)}(x, \cdot)$ together with the elliptic Harnack inequality   \eqref{e:2.5},
there exists $c_3=c_3(d, \phi, \Lambda_0)>0$ such that
\begin{equation}\label{e:2.8'}
G_{U_r(z_0)}(x, y)\leq c_3G_{U_r(z_0)}(x, y_0)\quad \mbox{for} \quad y\in U_r(z_0) \setminus D_{r/(2L)}(z_0) \quad \mbox{with} \quad \delta_{U_r(z_0)}(y)\leq r/(16L).
\end{equation}
On the other hand, note that $\delta_{U_r(z_0)}(y_0)\geq \kappa r/(2L).$
As $U_r(z_0)$  is above the graph of the local Lipschitz function $\Gamma_{z_0}$ in the coordinate system $CS_{z_0},$
  by  the Harnack chain argument and the uniform elliptic Harnack inequality \eqref{e:2.5}, there exists $c_4=c_4(d, \phi,\Lambda_0)$ such that
\begin{equation}\label{e:2.8}
G_{U_r(z_0)}(x, y)\leq c_4G_{U_r(z_0)}(x, y_0) \quad \mbox{for} \quad y\in U_r(z_0) \setminus D_{r/(2L)}(z_0) \quad \mbox{with} \quad \delta_{U_r(z_0)}(y)\geq r/(16L).
\end{equation}
By combining \eqref{e:2.8'} and \eqref{e:2.8},
\begin{equation}\label{e:2.9'}
G_{U_r(z_0)}(x, y)\leq c_5G_{U_r(z_0)}(x, y_0)\quad \mbox{for} \quad y\in U_r(z_0) \setminus D_{r/(2L)}(z_0),
\end{equation}
where $c_5:=c_3\vee c_4.$

For each $n\geq 2,$ let $D_n:=D\cap B(z_0, n).$
 By Ito's formula, for $n\geq 2$ and $x\in D_{r/(8L)}(z_0),$
\begin{equation}\label{e:2.9n'}\begin{aligned}
\P_x(X_{\tau_{U_r(z_0)}}\in D_n)&\leq \E_x \varphi^{(n)}_{z_0, r}(X_{\tau_{U_r(z_0)}})\\
&=\varphi^{(n)}_{z_0, r}(x)+\E_x \int_0^{\tau_{U_r(z_0)}} \L \varphi^{(n)}_{z_0, r}(X_s) \,ds\\
&=\int_{U_r(z_0)} G_{U_r(z_0)}(x, y)\L \varphi^{(n)}_{z_0, r}(y)\,dy\\
&= \int_{U_r(z_0)\setminus D_{r/(2L)}(z_0)} G_{U_r(z_0)}(x, y) \L \varphi^{(n)}_{z_0, r} (y)\,dy\\
&\quad+\int_{D_{r/(2L)}(z_0)} G_{U_r(z_0)}(x, y) \L \varphi^{(n)}_{z_0, r} (y)\,dy.
\end{aligned}\end{equation}
Note that $\varphi_{z_0, r}^{(n)}(y)=\varphi_{z_0, r}(y)=0$ for $y\in B(z_0, 3r/(4L)),$ by \eqref{e:2.17n}, we have
\begin{equation}\label{e:2.24}
\L \varphi^{(n)}_{z_0, r} (y)=\int_{\R^d \setminus B(z_0, 3r/(4L))}\dfrac{\varphi^{(n)}_{z_0, r}(z)}{|y-z|^d \phi(|y-z|)}\,dz \quad \mbox{for}\quad  y\in  D_{r/(2L)}(z_0).
\end{equation}
Then by \eqref{e:2.9n'} together with \eqref{e:2.7}, \eqref{e:2.9'}, \eqref{e:2.24} and the condition \eqref{e:1.1},
\begin{equation}\label{e:3.20n}\begin{aligned}
&\quad\P_x(X_{\tau_{U_r(z_0)}}\in D_n)\\
&\leq c_2c_5r^{-2}\int_{U_r(z_0)\setminus D_{r/(2L)}(z_0)} G_{U_r(z_0)}(x, y_0) \,dy\\
&\quad +\int_{D_{r/(2L)}(z_0)}G_{U_r(z_0)}(x, y) \left[\int_{\R^d \setminus B(z_0, 3r/(4L))} \dfrac{\varphi^{(n)}_{z_0, r}(z)}{|y-z|^d \phi(|y-z|)}\,dz\right]\,dy\\
&\leq c_6\dfrac{G_{U_r(z_0)}(x, y_0)}{r^{2-d}}\\
&\quad +c_6\E_x \tau_{U_r(z_0)} \cdot \int_{\R^d \setminus B(z_0, 3r/(4L))} \dfrac{1}{|z-z_0|^d \phi(|z-z_0|)}\,dz,
\end{aligned}\end{equation}
where $c_6=c_6(d, \phi, \Lambda_0).$
By the condition \eqref{e:1.1}, $\frac{\phi(r)}{\phi(|z-z_0|)}\leq c(r/|z-z_0|)^{\beta_2}1_{\{|z-z_0|\leq r\}}(z)+c(r/|z-z_0|)^{\beta_1}1_{\{|z-z_0|> r\}}(z).$
Thus there exists $c_7=c_7(d, \phi, \Lambda_0)>1$ such that
$$\begin{aligned}
&\int_{\R^d \setminus B(z_0, 3r/(4L))} \dfrac{\phi(r)}{|z-z_0|^d \phi(|z-z_0|)}\,dz\\
\leq &\int_{\R^d} [ cr^{\beta_2}|z-z_0|^{-(d+\beta_2)}1_{\{3r/(4L)<|z-z_0|\leq r\}}(z)+cr^{\beta_1}|z-z_0|^{-(d+\beta_1)}1_{\{|z-z_0|> r\}}(z) ]\,dz\\
\leq & c_7.
\end{aligned}$$
Hence by this inequality together with \eqref{e:3.20n},
\begin{equation}\label{e:2.9n}
\P_x(X_{\tau_{U_r(z_0)}}\in D_n)\leq c_6\dfrac{G_{U_r(z_0)}(x, y_0)}{r^{2-d}}+c_6c_7\dfrac{\E_x\tau_{U_r(z_0)}}{\phi(r)}.
\end{equation}
 Observe that $G_{U_r(z_0)}(\cdot, y_0)$ is regular harmonic in $D_{r/(4L)}(z_0).$ Then for $x\in D_{r/(8L)}(z_0),$
 \begin{equation}\label{e:2.10'}\begin{aligned}
 G_{U_r(z_0)}(x, y_0)&=\E_x [G_{U_r(z_0)}(X_{\tau_{D_{r/(4L)}(z_0)}}, y_0)]\\
 &\geq \E_x [G_{U_r(z_0)}(X_{\tau_{D_{r/(4L)}(z_0)}}, y_0); X_{\tau_{D_{r/(4L)}(z_0)}}\in B(y_0, \kappa r/(8L))]\\
 &\geq \inf_{y\in B(y_0, \kappa r/(8L))}G_{U_r(z_0)}(y, y_0)\P_x(X_{\tau_{D_{r/(4L)}(z_0)}}\in B(y_0, \kappa r/(8L)))\\
 &\geq \inf_{y\in B(y_0, \kappa r/(8L))}G_{B(y_0, \kappa r/(4L))}(y, y_0)\P_x(X_{\tau_{D_{r/(4L)}(z_0)}}\in B(y_0, \kappa r/(8L)))\\
&\geq c_8r^{2-d}\P_x(X_{\tau_{D_{r/(4L)}(z_0)}}\in B(y_0, \kappa r/(8L))),
\end{aligned}\end{equation}
where the last inequality is due to Lemma \ref{L:2.4'}.
By the condition \eqref{e:1.1} and the L\'evy system formula of $X,$
\begin{equation}\label{e:2.15}\begin{aligned}
\P_x(X_{\tau_{D_{r/(4L)}(z_0)}}\in B(y_0, \kappa r/(8L)))
&=\E_x \int_0^{\tau_{D_{r/(4L)}(z_0)}}\int_{B(y_0, \kappa r/(8L))} j(|X_s-y|) \,dy\,ds\\
&=\E_x \int_0^{\tau_{D_{r/(4L)}(z_0)}}\int_{B(y_0, \kappa r/(8L))} \dfrac{1}{|X_s-y|^d \phi(|X_s-y|)} \,dy\,ds\\
&\geq c_9\E_x \int_0^{\tau_{D_{r/(4L)}(z_0)}}\int_{B(y_0, \kappa r/(8L))} \dfrac{1}{r^d \phi(r)} \,dy\,ds\\\
&\geq c_{10}\dfrac{\E_x\tau_{D_{r/(4L)}(z_0)}}{\phi(r)},
\end{aligned}\end{equation}
where the third line holds as the distance between $D_{r/(4L)}(z_0)$ and $B(y_0, \kappa r/(8L))$ is less than $r.$
Hence by \eqref{e:2.10'} and \eqref{e:2.15},
\begin{equation}\label{e:2.11}
 \dfrac{G_{U_r(z_0)}(x, y_0)}{r^{2-d}}\geq c_8c_{10}\dfrac{\E_x\tau_{D_{r/(4L)}(z_0)}}{\phi(r)}.
 \end{equation}
 By \cite[Theorem 2.9]{CW2}, there exists $c_{11}=c_{11}(d, \phi)$ such that
\begin{equation}\label{e:2.17}
\E_x \tau_{U_r(z_0)}\leq \E_x \tau_{D_r(z_0)}\leq c_{11}\E_x\tau_{D_{r/(4L)}(z_0)}.
\end{equation}
By combining \eqref{e:2.9n}, \eqref{e:2.11} and \eqref{e:2.17}, there exists $c_{12}=c_{12}(d, \phi, \Lambda_0)$ such that for any $n\geq 2,$
$$\P_x(X_{\tau_{U_r(z_0)}}\in D_n)\leq c_{12}\dfrac{G_{U_r(z_0)}(x, y_0)}{r^{2-d}}.$$
By letting $n\rightarrow\infty,$ the desired result is obtained.

\qed

\begin{lem}\label{L:2.6}
Let $D$ be a  Lipschitz open set in $\R^d$ with characteristics $(R_0, \Lambda_0).$ For each $b\in (0, (2L)^{-1}),$
there exists $C=C(d, \phi, \kappa, b)>0$  such that for any $z_0\in \partial D$ and   $r\in (0, R_0),$
$$
\P_x \left( X_{\tau_{U_r(z_0)\setminus \overline{B(y_0, b\kappa r)}}}\in \overline{B(y_0, b\kappa r)} \right)
\geq C r^{d-2} G_{U_r(z_0)} (x, y_0) \quad  \hbox{for } x\in D\cap B(z_0, r/(8L)),
$$
where   $y_0 \in D\cap \partial B(z_0, r/(2L))$ with $\kappa r/(2L)< \delta_D(y_0)< r/(2L).$
\end{lem}

\pf The proof is similar to \cite[Lemma 3.10]{CW2}.
Let $y_0 \in D\cap \partial B(z_0, r/(2L))$ with $\kappa r/(2L)< \delta_D(y_0)< r/(2L).$
 Let $F:=U_r(z_0)\setminus \overline{B(y_0, b\kappa r)}.$ Since $G_{U_r(z_0)}(\cdot, y_0)$ is regular harmonic in $F,$ we have for $x\in D\cap B(z_0, r/(8L)),$
\begin{equation}\label{e:2.16}\begin{aligned}
G_{U_r(z_0)}(x, y_0)&=\E_x \left[ G_{U_r(z_0)}(X_{\tau_F}, y_0); X_{\tau_F}\in \overline{B(y_0, b\kappa r)} \right]\\
&=\E_x  \left[ G_{U_r(z_0)}(X_{\tau_F}, y_0); X_{\tau_F}\in B(y_0, b\kappa r/2) \right]\\
&\quad+\E_x  \left[ G_{U_r(z_0)}(X_{\tau_F}, y_0); X_{\tau_F}\in \overline{B(y_0, b\kappa r)}\setminus B(y_0, b\kappa r/2) \right]\\
&\leq \E_x\int_0^{\tau_F}\int_{B(y_0, b\kappa r/2)} G_{U_r(z_0)}(u, y_0) j(|X_s-u|)  \,du\,ds\\
&\quad+\sup_{u\in \overline{B(y_0, b\kappa r)}\setminus B(y_0, b\kappa r/2)}  G_{B(u, 4r)}(u, y_0)\cdot
\P_x \left( X_{\tau_F}\in \overline{B(y_0, b\kappa r)}\setminus B(y_0, b\kappa r/2) \right)\\
&\leq c_1 j(r) \E_x\tau_F\int_{B(y_0, b\kappa r/2)}G_{U_r(z_0)}(y_0, u)\,du\\
&\quad+c_2r^{2-d}\P_x \left( X_{\tau_F}\in \overline{B(y_0, b\kappa r)}\setminus B(y_0, b\kappa r/2) \right),
\end{aligned}\end{equation}
where $c_1=c_1(d, \phi, b, \kappa)$ and $c_2=c_2(d, \phi),$ the condition \eqref{e:1.1} and Lemma \ref{L:2.4'} are used in the last inequality.

 Note that by Lemma \ref{L:2.2}, there exists a positive constant $c_3=c_3(d, \phi)>1$ such that
\begin{equation}\label{e:3.42}
\int_{B(y_0, b\kappa r/2)}G_{U_r(z_0)}(y_0, u)\,du\leq \int_{B(y_0, 4r)}G_{B(y_0, 4r)}(y_0, u)\,du=\E_{y_0} \tau_{B(y_0, 4r)}\leq c_3r^2.
\end{equation}
On the other hand, by the condition \eqref{e:1.1},  there exists a positive constant $c_4=c_4(d, \phi, b, \kappa)>1$ such that for any $r\in (0, 1),$
\begin{equation}\label{e:3.43}\begin{aligned}
\P_x(X_{\tau_F}\in B(y_0, b\kappa r/2))
=\int_{B(y_0, b\kappa r/2)}\int_F G_F(x, y) j(|y-u|) \,dy\,du
\geq c_4  r^d j(r)  \E_x \left[ \tau_F \right].
\end{aligned}\end{equation}
Hence, by \eqref{e:3.42} and \eqref{e:3.43}, the first item on the right side of \eqref{e:2.16} satisfies
\begin{equation}\label{e:3.39'}\begin{aligned}
c_1 j(r) \E_x\tau_F\int_{B(y_0, b\kappa r/2)}G_{U_r(z_0)}(y_0, u)\,du
\leq c_1c_3 r^2j(r) \E_x\tau_F
\leq c_1c_3 c_4^{-1}r^{2-d}\P_x(X_{\tau_F}\in B(y_0, b\kappa r/2)).
\end{aligned}\end{equation}
Let $c_5:=(c_1c_3 c_4^{-1})\vee c_2.$ By \eqref{e:2.16} together with \eqref{e:3.39'}, we have
$$
G_{U_r(z_0)}(x, y_0) \leq 2c_5r^{2-d}\P_x(X_{\tau_{F}}\in \overline{B(y_0, b\kappa r)}).
$$
That is,
\begin{equation}\label{e:3.20}
\P_x \left( X_{\tau_{U_r(z_0)\setminus \overline{B(y_0, b\kappa r)}}}\in \overline{B(y_0, b\kappa r)} \right)
\geq (2c_5)^{-1} r^{d-2} G_{U_r(z_0)} (x, y_0) , \quad x\in D\cap B(z_0, r/(8L)).
\end{equation}
Hence, the proof is complete.
\qed

\begin{prp}\label{P:2.2}
Let $D$ be a Lipschitz open set in $\R^d$ with characteristics $(R_0, \Lambda_0).$

{\rm (i)} For each  $T>0,$  there exists $C_1=C_1(d, \phi, R_0, \Lambda_0, T)>1$  such that for any $z_0\in \partial D$ and $t\in (0, T),$
\begin{equation}\label{e:2.21}
\P_x(X_{\tau_{U_{R_0\wedge \sqrt t}(z_0)}}\in D)\leq C_1 \P_x(\tau_D>t), \quad  x\in D\cap B(z_0, (R_0\wedge\sqrt{t})/(8L)).
\end{equation}

{\rm(ii)} For each $T>0$ and $b\in (0, 1],$ there exists $C_2=C_2(d, \phi, b,  R_0, \Lambda_0, T)>1$  such that for any $z_0\in \partial D, t\in (0, T)$ and $s\in (0, t],$
\begin{equation}\label{e:2.22}
\P_x(\tau_{U_{b(R_0\wedge\sqrt{t})}(z_0)}>s)\leq  \P_x(\tau_D>s) \leq C_2\P_x(\tau_{U_{b(R_0\wedge\sqrt{t})}(z_0)}>s), \quad  x\in D\cap B(z_0, b(R_0\wedge\sqrt{t})/(8L)).
\end{equation}
\end{prp}

\pf  Let $t\in (0, T)$ and $t_0:=R_0\wedge \sqrt t.$
 Let $z_0\in \partial D$ and $b\in (0, 1].$ Let $x_0$ be a point of $D\cap \partial B(z_0, bt_0/(2L))$ with $\delta_D(x_0)\geq  b\kappa t_0/(2L).$ Note that $D\cap B(z_0, bt_0/L) \subset U_{bt_0}(z_0) \subset D\cap B(z_0, bt_0),$ thus $\delta_{U_{bt_0}(z_0)}(x_0)\geq  b\kappa t_0/(2L).$
Denote by $F:=U_{bt_0}(z_0)\setminus \overline{B(x_0, b\kappa t_0/(8L))}.$
Note that $t_0^2\leq t\leq (1\vee \frac{T}{R_0^2})t_0^2.$
By the strong Markov property of $X$, Lemmas \ref{L:2.7} and  \ref{L:2.6}, for  $x\in D\cap B(z_0, bt_0/(8L)),$
 \begin{equation}\label{e:2.10}\begin{aligned}
&\quad \P_x(\tau_{U_{bt_0}(z_0)}>t)\\
&\geq \E_x[ X_{\tau_F}\in \overline{B(x_0, b\kappa t_0/(8L))}; \P_{X_{\tau_F}} (\tau_{B(X_{\tau_F}, b\kappa t_0/(4L))}>t)]\\
&\geq \P_x[X_{\tau_F}\in \overline{B(x_0, b\kappa t_0/(8L))}] \inf_{y\in \overline{B(x_0, b\kappa t_0/(8L))}}\P_y(\tau_{B(y, b\kappa t_0/(4L))}>t)\\
&\geq c_1G_{U_{bt_0}(z_0)}(x, x_0)t_0^{d-2},
\end{aligned}\end{equation}
where $c_1=c_1(d, \phi, b, R_0, \Lambda_0,  T).$
By Proposition \ref{P:2.1}, there exists $c_2=c_2(d, \phi, b, \Lambda_0)$ such that
$$\P_x(X_{\tau_{U_{bt_0}(z_0)}}\in D)\leq c_2G_{U_{bt_0}(z_0)}(x, x_0)t_0^{d-2}, \quad  x\in D\cap B(z_0, bt_0/(8L)).$$
Hence by combining this inequality and \eqref{e:2.10},
\begin{equation}\label{e:2.38}
\P_x(X_{\tau_{U_{bt_0}(z_0)}}\in D)\leq c_2c_1^{-1}\P_x(\tau_{U_{bt_0}(z_0)}>t)\leq c_2c_1^{-1}\P_x(\tau_D>t), \quad  x\in D\cap B(z_0, bt_0/(8L)).
\end{equation}
Take $b=1$ in \eqref{e:2.38}, then \eqref{e:2.21} holds.

It follows from the first inequality in \eqref{e:2.38} that for $x\in D\cap B(z_0, bt_0/(8L))$ and $s\in (0, t],$
\begin{equation}\label{e:3.35}\begin{aligned}
\P_x(\tau_D>s)
&\leq \P_x(\tau_{U_{bt_0}(z_0)}>s)+\P_x(X_{\tau_{U_{bt_0}(z_0)}}\in D)\\
&\leq \P_x(\tau_{U_{bt_0}(z_0)}>s)+c_2c_1^{-1}\P_x(\tau_{U_{bt_0}(z_0)}>t)\\
&\leq (1+c_2c_1^{-1})\P_x(\tau_{U_{bt_0}(z_0)}>s).
\end{aligned}\end{equation}
It is easy to see that for $z_0\in \partial D$ and $s\in (0, t),$
\begin{equation}\label{e:3.36}
\P_x(\tau_{U_{bt_0}(z_0)}>s)\leq  \P_x(\tau_D>s).
\end{equation}
Hence by \eqref{e:3.35} and \eqref{e:3.36}, \eqref{e:2.22} holds. Thus the proof is complete.
\qed

\section{Parabolic Carleson estimate and survival probability estimates}

In this section, we shall establish  parabolic Carleson estimate for a class of nonnegative parabolic functions with respect to $Z_t=(V_0-t, X_t)$ in cylindrical set $(0, 4r^2)\times U$  that vanish on $(0, 4r^2)\times U^c$ for $r\in (0, T)$ and a small Lipschitz subset $U$ of $D$ near the boundary, where $T$ is a positive constant.  Furthermore we prove  the  doubling property of the survival probability for $X$ in  Lipschitz open sets.

Let
  $Z=\{Z_t:=(V_0-t, X_t), t\geq 0\}$
 be the
time-space process of $X$. Let $D$ be an open set.
Recall that
a non-negative real valued Borel  measurable function
$h(t,x)$ on $[0, \infty)\times \R^d$ is {\it parabolic}
(or {\it caloric})
on $Q=(a,b)\times D \subset [0,\infty)\times\R^d$
if
for every relatively compact open subset $Q_1$ of $Q$,
$$
h(t, x)=\E_{(t,x)} [h (Z_{\tau_{Q_1}})]
$$
for every $(t, x)\in Q_1$,
where
$\tau_{Q_1}=\inf\{s> 0: \, Z_s\notin Q_1\}$.

\begin{lem}\label{L:3.1}
Let $D$ be a Lipschitz open set.
For each $y\in D,$ $(t, x)\mapsto p_D(t, x, y)$ is  parabolic with respect to $Z=\{Z_t:=(V_0-t, X_t), t\geq 0\}$ on $(0, \infty)\times D$ and vanishes continuously  on
$(0, \infty)\times D^c.$
\end{lem}

\pf By the strong Markov property of $X,$ for every relatively compact open subset $U\subset D,$
\begin{equation}\label{e:3.1}
p_D(t, x, y)=p_U(t, x, y)+\E_x [p_D(t-\tau_U, X_{\tau_U}, y); \tau_U<t], \quad x, y\in D, t>0.
\end{equation}
Let $U$ be a relatively compact open subset in $D.$ Define $Q_1:=(a, b)\times U\subset (0, \infty)\times D$ and $\tau_1:=\inf\{s>0: Z_s\notin Q_1\}.$
It is easy to see that $\tau_1=(t-a) \wedge \tau_U$ for starting point $Z_0=(t, x)\in Q_1.$

Fix $y\in D.$ Let $h(t, x):=p_D(t, x, y).$ Then for $(t, x)\in Q_1=(a, b)\times U,$
\begin{equation}\label{e:3.2}\begin{aligned}
&\E_{(t, x)} h(Z_{\tau_1})\\
=&\E_x p_D(t-\tau_1, X_{\tau_1}, y)\\
=&\E_x [p_D(t-\tau_1, X_{\tau_1}, y); \tau_U<t-a]+\E_{(t, x)} [p_D(t-\tau_1, X_{\tau_1}, y); \tau_U\geq t-a]\\
=&\E_x [p_D(t-\tau_U, X_{\tau_U}, y); \tau_U<t-a]+\E_x[p_D(a, X_{t-a}, y); \tau_U\geq t-a]\\
=&\E_x [p_D(t-\tau_U, X_{\tau_U}, y); \tau_U<t-a]+\E_x[p_U(a, X_{t-a}, y); \tau_U\geq t-a]\\
&+\E_x[\E_{X_{t-a}}p_D(a-\tilde \tau_U, X_{\tilde \tau_U}, y); \tau_U\geq t-a],
\end{aligned}\end{equation}
where  \eqref{e:3.1} is used in the last equality and $\tilde \tau_U$ is the first exit time of $X$ from $U$ starting from $X_{t-a}$.

Let $X^U$ be the subprocess of $X$ killed upon exiting $U.$
Note that
\begin{equation}\label{e:3.3}\begin{aligned}
&\E_x[p_U(a, X_{t-a}, y); \tau_U\geq t-a]= \E_x[p_U(a, X^U_{t-a}, y)]\\
=&\int_U p_U(t-a, x, z) p_U(a, z, y)\,dz=p_U(t, x, y).
\end{aligned}\end{equation}
Let $\theta_t$ be the time shift operator of $X.$ By the Markov property of $X,$ we have
\begin{equation}\label{e:3.4}\begin{aligned}
&\E_x[\E_{X_{t-a}}p_D(a-\tilde\tau_U, X_{\tilde\tau_U}, y); \tau_U\geq t-a]\\
=&\E_x [\E_x(\theta_{t-a}\circ p_D(a-\tilde\tau_U, X_{\tilde\tau_U}, y)|\mathcal{F}_{t-a}); \tau_U\geq t-a]\\
=&\E_x [\E_x(p_D(t-\tau_U, X_{\tau_U}, y); \tau_U\geq t-a|\mathcal{F}_{t-a})]\\
=&\E_x [p_D(t-\tau_U, X_{\tau_U}, y); \tau_U\geq t-a].
\end{aligned}\end{equation}
Then by combining \eqref{e:3.2} with \eqref{e:3.3} and \eqref{e:3.4},
\begin{equation}\label{e:3.5'}\begin{aligned}
&\E_{(t, x)} h(Z_{\tau_1})\\
=& \E_x [p_D(t-\tau_U, X_{\tau_U}, y); \tau_U<t-a]+p_U(t, x, y)+\E_x [p_D(t-\tau_U, X_{\tau_U}, y); \tau_U\geq t-a]\\
=& p_U(t, x, y)+\E_x p_D(t-\tau_U, X_{\tau_U}, y)\\
=& p_D(t, x, y)=h(t, x).
\end{aligned}\end{equation}
Hence $(t, x)\mapsto p_D(t, x, y)$ is  parabolic with respect to $Z$ on $(0, \infty)\times D.$

On the other hand, since $D$ is a Lipschitz open set, it satisfies the exterior cone condition.
By Theorem \ref{T1} and  an argument similar to that of \cite[Proposition 1.13, page 89]{Ba},  each boundary point of $\partial D$ is regular for $D^c.$
 Then by an argument similar to that of \cite[Theorem 2.4]{ChZ}, we have
\begin{equation}\label{e:3.6'}
\lim_{x\rightarrow z\in \partial D}p_D(t, x, y)=0.
\end{equation}

The proof is complete by \eqref{e:3.5'} and \eqref{e:3.6'}.

\qed

\smallskip

Let $D$ be a Lipschitz open set with characteristics $(R_0, \Lambda_0).$
Recall that for each $z\in\partial D$ and $r\in (0, R_0),$  $U_r(z)\subset D$  is a Lipschitz connected open set    with characteristics $(r/L, L\Lambda_0)$ such that $D\cap B(z, r/L)\subset U_r(z)\subset D\cap B(z, r)$, where $L=L(\Lambda_0, d)>1.$

\begin{lem}\label{L:3.2}
Let $D$ be a Lipschitz open set with characteristics $(R_0, \Lambda_0).$
For each $T>0,$
there exist $C=C(d, R_0, \Lambda_0, \phi, T)>0$ and $\beta=\beta(d, \Lambda_0, \phi)>0$ such that for any $z_0\in \partial D, r\in (0, T)$ and any  nonnegative function $h(t, x)$  which is parabolic with respect to $Z_t=(V_0-t, X_t)$ in $Q=(0, 4r^2)\times U_{r\wedge R_0}(z_0),$
\begin{equation}\label{e:3.8'}
h(s, x)\leq C\left(\frac{r\wedge R_0}{\delta_{U_{r\wedge R_0}(z_0)}(x)}\right)^{\beta} h(2r^2, x_0)
\end{equation}
for $s\in (r^2/4, r^2]$ and $x\in  U_{r\wedge R_0}(z_0)$ with $ \delta_{U_{r\wedge R_0}(z_0)}(x)  < (r\wedge R_0)/(16L),$
where $x_0$ is a point of $D\cap \partial B(z_0, (r\wedge R_0)/(2L))$ with $\delta_D(x_0)>\kappa (r\wedge R_0)/(2L).$
\end{lem}

\pf Let $r\in (0, T)$  and $z_0\in\partial D.$
For the simplicity of notation, denote by $U:=U_{r\wedge R_0}(z_0).$
Recall that  $U$ is a Lipschitz domain with characteristics $((r\wedge R_0)/L, L\Lambda_0).$
Fix $x\in  U$ with $\delta_U(x)< (r\wedge R_0)/(16L).$
 Let $z_x$ be a point of $\partial U$ such that $|x-z_x|=\delta_{U}(x).$
 Denote by $CS_{z_x}^U$  the local  coordinate system of $U$ with its origin at $z_x.$
Let $x=(x^1, \cdots, x^{d-1}, x^d)$ be the coordinate of $x$ in $CS_{z_x}^U.$
Define $y_1:=x=(x^1, \cdots, x^{d-1}, x^d)$ and $y_i:=(x^1, \cdots, x^{d-1}, y^d_i), i\geq 2$ in $CS_{z_x}^U$  such that
 $b_i:=|y_{i+1}-y_i|=y_i^d\frac{1}{4\sqrt{1+L^2\Lambda_0^2}}, i\geq 1.$ Let $n$ be the first integer such that $y_n^d> (r\wedge R_0)/(2L).$
Since $U$ is a Lipschitz domain with characteristics $((r\wedge R_0)/L, L\Lambda_0),$
 it is easy to see that $b_i\leq \frac{1}{4}{\rm dist} (y_i, \partial U)$ for $i=1, \cdots, n.$

 Let $\gamma\in (0, 1)$ be a constant which will be determined later.
For each $s\in (r^2/4, r^2],$   define $s_1=s$ and $s_{i+1}=s_i+\gamma b_i^2$ for $i\geq 1.$
Note that $y_{i+1}\in B(y_i, \frac{1}{4}{\rm dist} (y_i, \partial U))$ for $i\geq 1$ and $R_0\leq 1.$
Let $h(t, x)$  be a nonnegative parabolic function with respect to $Z_t=(V_0-t, X_t)$ in $Q=(0, 4r^2)\times U.$
By the parabolic Harnack inequality  in Proposition \ref{P:3.1}, there exists $c_1=c_1(d, \phi,  \gamma)>1$ such that for $i=1, \cdots, n-1,$
$$h(s_i, y_{i})\leq c_1h(s_{i+1}, y_{i+1}).$$
Hence,
\begin{equation}\label{e:3.9'}
h(s, x)=h(s_1, y_1)\leq c_1^n h(s_n, y_n).
\end{equation}
Let $a:=\frac{1}{4\sqrt{1+L^2\Lambda_0^2}}.$
Since $y_{i+1}^d=y_i^d (1+a)=x^d (1+a)^i$ for $i\geq 1,$
we have $b_i=y_{i+1}^d-y_i^d=ay_i^d=x^d a(1+a)^{i-1}.$
Thus
$$\begin{aligned}
s_n&=s+\sum_{i=1}^n \gamma^2 b_i^2=s+(x^d)^2a^2\sum_{i=1}^n \gamma^2 (1+a)^{2(i-1)}\\
&=s+(x^d)^2\gamma^2a^2(1+a)^{-2} \dfrac{(1+a)^{2n}-1}{(1+a)^2-1}\\
&\leq  s+\dfrac{(y^d_n)^2}{(1+a)^2-1}\gamma^2a^2.
\end{aligned}$$
Note that $n$ is the first integer such that $y_n^d> (R_0\wedge r)/(2L).$ Then $y^d_n \leq  (R_0\wedge r)(1+a)/(2L).$ We have
$$\begin{aligned}
s_n
&\leq r^2+\dfrac{((R_0\wedge r)/(2L))^2(1+a)^2}{(1+a)^2-1}\gamma^2a^2\leq r^2+\dfrac{( r/2)^2(1+a)^2}{(1+a)^2-1}\gamma^2a^2.
\end{aligned}$$
We take $\gamma=\gamma(d, \Lambda_0)$ small enough such that $s_n\leq 3r^2/2.$
Since $y_n^d=x^d(1+a)^{n-1}$ and $y^d_n \leq  (r\wedge R_0)(1+a)/(2L),$ it follows that $n\leq \log(\frac{(r\wedge R_0)}{x^d})/\log(1+a)+2.$
Thus by \eqref{e:3.9'}, there exist $\beta=\beta(d,  \Lambda_0, \phi)>0$ and $c_2=c_2(d, \Lambda_0, \phi)>0$ such that
\begin{equation}\label{e:3.9n}
h(s, x)\leq c_1^n h(s_n, y_n)\leq c_2 ((r\wedge R_0)/x^d)^{\beta}  h(s_n, y_n).
\end{equation}
Let $x_0$ be a point of $D\cap \partial B(z_0, (r\wedge R_0)/(2L))\subset U$ with $\delta_U(x_0)=\delta_D(x_0)>\kappa (r\wedge R_0)/(2L).$
Note that $s_n\in (s, 3r^2/2)\subset (r^2/4, 3r^2/2)$ and $(r\wedge R_0)\leq r\leq (1\vee \frac{T}{R_0}) (r\wedge R_0)$ for $r\in (0, T).$
By the chain argument and the parabolic Harnack inequality in Proposition \ref{P:3.1}, there exists  $c_3=c_3(d,  \phi, R_0, \Lambda_0, T)>0$  such that for $r\in (0, T),$
\begin{equation}\label{e:3.10n}
h(s_n, y_n)\leq c_3 h(2r^2, x_0).
\end{equation}
Consequently, \eqref{e:3.8'} holds by \eqref{e:3.9n} and \eqref{e:3.10n}.

\qed

\begin{prp}[Parabolic Carleson estimate]\label{P:3.2}
Let $D$ be a Lipschitz open set with characteristics $(R_0, \Lambda_0).$
For each $T>0,$
there exists $C=C(d, R_0, \Lambda_0, \phi, T)>0$ such that for any $r\in (0, T), z_0\in \partial D$ and any  nonnegative function $h(t, x)$  which is parabolic with respect to $Z_t=(V_0-t, X_t)$ in $Q=(0, 4r^2)\times U_{r\wedge R_0}(z_0)$ and vanishes continuously  on $(0, 4r^2)\times U_{r\wedge R_0}(z_0)^c,$
\begin{equation}\label{e:3.5}
h(r^2, x)\leq C h(2r^2, x_0) \quad \mbox{for} \quad \: x\in  U_{r\wedge R_0}(z_0) \:\: \mbox{with} \:\: \delta_{U_{r\wedge R_0}(z_0)}(x)< (r\wedge R_0)/(16L),
\end{equation}
where $x_0$ is a point of $D\cap \partial B(z_0, (r\wedge R_0)/(2L))$ with $\delta_{U_{r\wedge R_0}(z_0)}(x_0)>\kappa (r\wedge R_0)/(2L).$
\end{prp}

\pf Let $r\in (0, T)$ and $z_0\in\partial D.$
For the simplicity of notation, denote by $U:=U_{r\wedge R_0}(z_0).$
Let $r_1:=(r\wedge R_0)/L$ and $\Lambda_1:=\Lambda_0 L.$
Note that $U$ is a Lipschitz domain with characteristics $(r_1, \Lambda_1).$
 Let $x_0$ be a point of $D\cap \partial B(z_0, r_1/2)$ with $\delta_{U}(x_0)>\kappa r_1/2.$

 Let $h(t, x)$  be a parabolic function with respect to $Z_t=(V_0-t, X_t)$ in $Q=(0, 4r^2)\times U$ and vanishes on $(0, 4r^2)\times U^c.$
 Let $K>4$ be a constant which will be determined later.
 For the simplicity of notation, we denote $U_{r_1/(16K)}:=\{x\in U: \delta_U(x)<r_1/(16K)\}.$
 Let $x\in U_{r_1/(16K)}$  and
\begin{equation}\label{e:3.9}
B_0:=U\cap B(x, 4\delta_U(x)),\,\quad  B_1:=U\cap B(x,  r_1/8),\qquad B_2:=B(x_0, \delta_U(x_0)/4).
\end{equation}
Set
\begin{equation}\label{e:3.10}
 \tau_0:=\tau_{B_0}\wedge K\delta_U(x)^2\, , \qquad  \tau_2:=\tau_{B_2}\wedge \dfrac{7}{4}r^2.
 \end{equation}
Then for $x\in U_{r_1/(16K)}$ and $t\in (r^2/2, r^2),$
$$\begin{aligned}
h(t, x)&=\E_{(t, x)} h(Z_{\tau_0})=\E_x h(t-\tau_0, X_{\tau_0})\\
&=\E_x [h(t-\tau_0, X_{\tau_0}); X_{\tau_0}\in B_1]
+\E_x [h(t-\tau_0, X_{\tau_0}); X_{\tau_0}\in B_1^c\cap U].
\end{aligned}$$

We first show that there exists $c_0=c_0(d, R_0, \phi, T)>0$ such that for $x\in U_{r_1/(16K)}$ and $t\in (r^2/2, r^2),$
\begin{equation}\label{e:3.7}
\E_x [h(t-\tau_0, X_{\tau_0}); X_{\tau_0}\in B_1^c\cap U]\leq c_0h(2r^2, x_0) .
\end{equation}
Note that the distance between $B_0$ and $B_1^c$ is larger than $r_1/16.$
By the L\'evy system formula \eqref{e:2.2} and the condition \eqref{e:1.1},
\begin{equation}\label{e:3.15}\begin{aligned}
&\E_x\big[h(t-\tau_0, X_{\tau_0}); X_{\tau_0}\in B_1^c\cap U\big]
\\
=& \E_x \int_0^{\tau_0}\int_{B_1^c\cap U} h(t-s, y)j(|X_s-y|)\,dy\,ds\\
\leq & c_1\E_x \int_0^{\tau_0}\int_{B_1^c\cap U} h(t-s, y) j(r_1)\,dy\,ds\\
\leq  & c_1\int_0^{K\delta_U^2(x)}\int_{B_1^c\cap U} h(t-s, y)j(r_1)\,dy\,ds\\
= & c_1 \int_0^{K\delta_U^2(x)}
\left(\int_{B_1^c\cap U\cap \{y: |y-x_0|>2
\delta_U(x_0)/3\}} j(r_1)h(t-s, y)\, dy\right. \\
  &\quad \quad \quad \quad \quad +\left.\int_{\{y: |y-x_0|\le
2\delta_U(x_0)/3\}}j(r_1)h(t-s, y)\, dy\right) \,ds\\
=:& c_1 (I_1+I_2).
\end{aligned}\end{equation}

For the first item $I_1,$  as $K\delta_U^2(x)\leq r_1^2/16\leq r^2/16$ for $x\in U_{r_1 /(16K)},$
we have for $t\in (r^2/2, r^2),$
\begin{equation}\label{e:3.11'}\begin{aligned}
I_1 &=  \int_{t-K\delta_U^2(x)}^t\int_{U\cap B_1^c\cap \{|y-x_0|>2
\delta_U(x_0)/3\}} j(r_1)h(s, y)\, dy\, ds\\
 &\leq  \int_{r^2/4}^{r^2}\int_{U\cap B_1^c\cap \{|y-x_0|>2
\delta_U(x_0)/3\}} j(r_1)h(s, y)\, dy\, ds.
\end{aligned}\end{equation}
On the other hand, by the L\'evy system formula \eqref{e:2.2} and the condition \eqref{e:1.1}, there exist $c_k=c_k(d, R_0, \phi, T)>0, k=2, 3$ such that
\begin{equation}\label{e:3.16'}\begin{aligned}
h(2r^2, x_0)&=\E_{x_0} h(2r^2-\tau_2, X_{\tau_2})\\
&\geq \E_{x_0} [h(2r^2-\tau_2, X_{\tau_2}); |X_{\tau_2}-x_0|>2
\delta_U(x_0)/3, X_{\tau_2}\in U\setminus B_1]\\
&=\E_{x_0} \int_0^{\tau_2}\int_{(U\setminus B_1)\cap \{y: |y-x_0|>2\delta_U(x_0)/3\}} h(2r^2-s, y) j(|X_s-y|)\,dy\,ds\\
&\geq c_2\E_{x_0} \int_0^{\tau_2}\int_{(U\setminus B_1)\cap \{y: |y-x_0|>2\delta_U(x_0)/3\}} h(2r^2-s, y) j(r_1)\,dy\,ds\\
&\geq c_2\int_0^{7r^2/4}\int_{B_1^c\cap U\cap \{y: |y-x_0|>2\delta_U(x_0)/3\}} h(2r^2-s, y)j(r_1)\P_{x_0}(\tau_{B_2}>s)\,dy\,ds\\
&\geq c_2\int_0^{7r^2/4}\int_{B_1^c\cap U\cap\{y: |y-x_0|>2\delta_U(x_0)/3\}} h(2r^2-s, y)j(r_1)\P_{x_0}(\tau_{B_2}>2r^2)\,dy\,ds\\
&\geq c_3\int_0^{7r^2/4}\int_{B_1^c\cap U\cap\{y: |y-x_0|>2\delta_U(x_0)/3\}} h(2r^2-s, y)j(r_1)\,dy\,ds\\
&= c_3\int_{r^2/4}^{2r^2}\int_{B_1^c\cap U\cap \{y: |y-x_0|>2\delta_U(x_0)/3\}} h(s, y)j(r_1)\,dy\,ds,
\end{aligned}\end{equation}
where the second line  from the bottom  is due to  Lemma \ref{L:2.7} and $r_1\leq r\leq L(1\vee \frac{T}{R_0})r_1.$
Hence it follows from \eqref{e:3.11'} and \eqref{e:3.16'} that
\begin{equation}\label{e:3.14}
I_1\leq c_3^{-1}h(2r^2, x_0).
\end{equation}

For the second item $I_2,$
note that $t-s\in (r^2/4, r^2)$ for $t\in (r^2/2, r^2)$ and $s\in (0, K\delta_U^2(x))\subset (0, r^2/16).$
By the parabolic Harnack inequality in Proposition \ref{P:3.1} and \eqref{e:1.4}, we have  for $t\in (r^2/2, r^2),$
\begin{equation}\label{e:3.16}\begin{aligned}
I_2 &= \int_0^{K\delta_U^2(x)}\int_{|y-x_0|\le 2\delta_U(x_0)/3}
 h(t-s, y)j(r_1)\, dy\,ds\\
&\leq c_4h(2r^2, x_0) \int_0^{K\delta_U^2(x)}\int_{|y-x_0|\le 2\delta_U(x_0)/3}
r_1^{-(d+\beta_2)}\, dy\,ds\\
&\leq c_5h(2r^2, x_0)K\delta_U^2(x)r_1^{-\beta_2}\\
&\leq c_5h(2r^2, x_0)K(r_1/16K)^2r_1^{-\beta_2}\\
&=c_5h(2r^2, x_0)\frac{r_1^{2-\beta_2}}{16^2 K}\leq c_5h(2r^2, x_0),
\end{aligned}\end{equation}
where $c_k=c_k(d, R_0, \phi, T)>0, k=4, 5.$
Consequently,  \eqref{e:3.7} holds by \eqref{e:3.15}, \eqref{e:3.14} and \eqref{e:3.16}.

\medskip

 We now use the contradiction method to prove that for each $T>0,$
there exists $C=C(d, R_0, \Lambda_0, \phi, T)>0$ such that for any $r\in (0, T), z_0\in \partial D$ and  $(t, x)\in (r^2/2, r^2]\times U_{r_1/(16K)},$
\begin{equation}\label{e:3.18'}
h(t, x)\leq Ch(2r^2, x_0).
\end{equation}
  Since $U$ is a Lipschitz domain with characteristics $(r_1, \Lambda_1)$, by an argument similar to that of \eqref{e:2.3}, there exists $\theta_1=\theta_1(d, \Lambda_1, \phi)\in (0, 1)$  such that for any $x\in U_{r_1/(16K)},$
\begin{equation}\label{e:3.7'}
\P_x(X_{\tau_{U\cap B(x, 4\delta_U(x))}}\in U)\leq \theta_1.
\end{equation}
By Lemma \ref{L:2.2}, we  choose $K=K(d, \phi, \theta_1)>4$ large enough such that  for any $x\in U_{r_1/(16K)},$
\begin{equation}\label{e:3.8}
\P_x (\tau_{B(x, 4\delta_U(x))}>K\delta_U(x)^2)\leq \dfrac{\E_x \tau_{B(x, 4\delta_U(x))}}{K\delta_U(x)^2}\leq (1-\theta_1)/2.
\end{equation}

Without loss of generality, we assume
 $h(2r^2, x_0)=1$.
Let $N>1$ be a large enough constant which will be determined later.
Suppose there exists $(t_1, x_1)\in (r^2/2, r^2]\times U_{r_1/(16K)}$  with $h(t_1, x_1)\geq N.$
Define $\tau_0$ and $B_k, k=0, 1$ as in \eqref{e:3.9} and \eqref{e:3.10} but with $x_1$ in place of $x.$
By \eqref{e:3.7}, we have
\begin{equation}\label{e:3.18}\begin{aligned}
h(t_1, x_1)&=\E_{x_1} h(t_1-\tau_0, X_{\tau_0})\\
&=\E_{x_1} [h(t_1-\tau_0, X_{\tau_0}); X_{\tau_0}\in B_1]
+\E_{x_1} [h(t_1-\tau_0, X_{\tau_0}); X_{\tau_0}\in B_1^c\cap U]\\
&\leq \sup_{(t, y)\in (t_1-K\delta_U(x_1)^2, t_1)\times B_1} h(t, y)\cdot \P_{x_1}(X_{\tau_0}\in B_1)+c_0 .
\end{aligned}\end{equation}
It follows from \eqref{e:3.7'} and \eqref{e:3.8} that
\begin{equation}\label{e:3.19}\begin{aligned}
&\quad\P_{x_1}(X_{\tau_0}\in B_1)\leq \P_{x_1}(X_{\tau_0}\in U)\\
&\leq \P_{x_1}(X_{\tau_{U\cap B(x_1, 4\delta_U(x_1))}}\in U; \tau_{U\cap B(x_1, 4\delta_U(x_1))}\leq K\delta_U(x_1)^2)+ \P_{x_1} (\tau_{U\cap B(x_1, 4\delta_U(x_1))}>K\delta_U(x_1)^2)\\
&\leq \P_{x_1}(X_{\tau_{U\cap B(x_1, 4\delta_U(x_1))}}\in U)+ \P_{x_1} (\tau_{B(x_1, 4\delta_U(x_1))}>K\delta_U(x_1)^2)\\
&\leq (1+\theta_1)/2.
\end{aligned}\end{equation}
Thus by combining \eqref{e:3.18} and \eqref{e:3.19},
\begin{equation}\label{e:3.25}
 h(t_1, x_1)\leq \dfrac{1+\theta_1}{2}\sup_{(t, y)\in (t_1-K\delta_U(x_1)^2, t_1)\times B_1}  h(t, y)+c_0.
\end{equation}
Note that $\frac{1+\theta_1}{2}\in (0, 1).$
Let $A>1$ be a constant such that $A(1+\theta_1)/2\in (0, 1).$
We take $N>1$  such that $N>c_0\frac{A}{A-1}.$
Hence by \eqref{e:3.25},
\begin{equation}\label{e:3.26n}
\dfrac{1+\theta_1}{2}\sup_{(t, y)\in (t_1-K\delta_U(x_1)^2, t_1)\times B_1}  h(t, y) \ge
h(t_1, x_1)-c_0 > \frac{h(t_1, x_1)}{A}.
\end{equation}
In the last inequality of the display above we used the assumption that
 $h(t_1, x_1)>N>\frac{A}{A-1}c_0$.
Let $\rho:=A(1+\theta_1)/2\in (0, 1).$ It follows from \eqref{e:3.26n} that there exist $x_2\in B_1=U\cap B(x_1, r_1/8)$ and $t_2\in (t_1-K\delta_U(x_1)^2, t_1)$ such that
$$
h(t_2, x_2)\geq h(t_1, x_1)/\rho\geq N/\rho.
$$
On the other hand, it follows from Lemma \ref{L:3.2} that  there exist  $c_6=c_6(d, \phi, R_0, \Lambda_0, T)$ and $\beta=\beta(d, \phi, \Lambda_0)$ such that for $r\in (0, T),$
$$h(t_1, x_1)\leq c_6(r_1/\delta_U(x_1))^{\beta} h(2r^2, x_0)= c_6(r_1/\delta_U(x_1))^{\beta}.$$
Hence,
$$\delta_U(x_1)\leq c_6^{1/\beta}r_1h(t_1, x_1)^{-1/\beta}\leq c_6^{1/\beta}r_1N^{-1/\beta}.$$
Thus
$$|t_1-t_2|\leq K\delta^2_U(x_1)\leq Kc_6^{2/\beta}r_1^2N^{-2/\beta}.$$

By induction, there exist a sequence $\{(t_i, x_i), i\geq 1\}$ with $t_{i+1}\in (t_i -K\delta_{U}(x_i)^2, t_i)$ and $x_{i+1}\in U\cap B(x_i, r_1/8)$ for each $i\geq 1$ such that
\begin{equation}\label{e:3.23'}
h(t_i, x_i)\geq N/\rho^{i-1}, \quad \delta_U(x_i)\leq c_6^{1/\beta}r_1h(t_i, x_i)^{-1/\beta}\leq c_6^{1/\beta}r_1(N/\rho^{i-1})^{-1/\beta}.
\end{equation}
$$|t_{i+1}-t_i|\leq K\delta_U^2(x_i)\leq Kc_6^{2/\beta}r_1^2(\rho^{i-1}/N)^{2/\beta}.$$
Hence, for each $i\geq 2,$
$$|t_i-t_1|\leq \sum_{j=2}^i |t_j-t_{j-1}|\leq \sum_{j=2}^i Kc_6^{2/\beta}r_1^2(\rho^{j-2}/N)^{2/\beta}\leq r^2Kc_6^{2/\beta}N^{-2/\beta}\sum_{j=0}^\infty \rho^{2j/\beta}.$$
By taking $N$ large enough such that
$Kc_6^{2/\beta}N^{-2/\beta}\sum_{j=0}^\infty \rho^{2j/\beta}<1/4.$
Note that $t_1\in (r^2/2, r^2].$
Thus  $t_i\in (r^2/4, r^2)$ for $i\geq 1.$
Observe that $h(t_i, x_i)\geq N/\rho^{i-1}\rightarrow\infty$ as $i\rightarrow\infty$ and $x_i\rightarrow\partial U$ as $i\rightarrow\infty$ by \eqref{e:3.23'}. This contradicts with the fact that $h$ vanishes continuously  on $(0, 4r^2)\times U^c$.
Hence this contradiction indicates that \eqref{e:3.18'} holds.
Finally, for $x\in  U$ with $r_1/(16K)<\delta_{U}(x)<r_1/16$, the conclusion \eqref{e:3.5} can be obtained by the parabolic Harnack inequality in Proposition \ref{P:3.1}.
Thus the proof  is complete.

\qed

\begin{prp}\label{P:3.3}
Let $D$ be a Lipschitz open set with characteristics $(R_0, \Lambda_0).$
For each $T>0,$ there exists $C=C(d, \phi, R_0, \Lambda_0, T)$ such that for $z\in \partial D, t\in (0, T)$ and $x, y\in  U_{\sqrt{t}\wedge R_0}(z)$  with $\delta_{U_{\sqrt{t}\wedge R_0}(z)}(y)\leq  (\sqrt{t}\wedge R_0)/(16L),$
$$p_{U_{\sqrt{t}\wedge R_0}(z)}(t, x, y)\leq C p_{U_{\sqrt{t}\wedge R_0}(z)}(2t, x, y_0),$$
where $y_0$ is a point of $D\cap \partial B(z, (\sqrt{t}\wedge R_0)/(2L))$ with $\delta_D(y_0)\geq \kappa (\sqrt{t}\wedge R_0)/(2L).$
\end{prp}

\pf For the simplicity of notation, we denote $U_{\sqrt{t}\wedge R_0}(z)$ by $U.$
By Lemma \ref{L:3.1} and the symmetry of $p_U(t, x, y)$ in $x$ and $y,$ for each $x\in U,$ $(t, y)\mapsto p_U(t, x, y)$ is parabolic with respect to $Z=\{Z_t:=(V_0-t, X_t), t\geq 0\}$ in $(0, \infty)\times U$ and vanishes on $(0, \infty)\times U^c.$
Hence the desired conclusion is obtained by Proposition \ref{P:3.2} with $\sqrt t$ in place of $r$ there.

\qed

\begin{prp}\label{P:3.5}
Let $D$ be a Lipschitz open set with characteristics $(R_0, \Lambda_0).$
For each $T>0,$ there exists $C=C(d, \phi, R_0, \Lambda_0, T)$ such that for $z_0\in \partial D$ and $t\in (0, T),$
$$\P_x(\tau_D>t) \leq C p_{U_{\sqrt t\wedge R_0}(z_0)}(2t, x, y_0)(\sqrt t \wedge R_0)^{d},  \quad  x\in D\cap B(z_0, (R_0\wedge\sqrt{t})/(8L)),$$
where $y_0$ is a point of $D\cap \partial B(z_0, (\sqrt{t}\wedge R_0)/(2L))$ with $\delta_D(y_0)\geq \kappa (\sqrt{t}\wedge R_0)/(2L).$
\end{prp}

 \pf Let $T>0$ and $t\in (0, T).$ Let $z_0\in \partial D.$
 For the simplicity of notation, we denote by $U_{\sqrt t\wedge R_0}:=U_{\sqrt t\wedge R_0}(z_0).$
By Proposition \ref{P:2.2} (ii), there exists $c_1=c_1(d, R_0, \Lambda_0, \phi, T)>0$ such that for $t\in (0, T),$
\begin{equation}\label{e:3.11}
\P_x(\tau_D>t)\leq c_1\P_x(\tau_{U_{\sqrt t\wedge R_0}}>t),  \quad  x\in D\cap B(z_0, (R_0\wedge\sqrt{t})/(8L)).
\end{equation}
Let $y_0$ be a point of $D\cap \partial B(z_0, (\sqrt t\wedge R_0)/(2L))$ with $\delta_D(y_0)\geq \kappa (\sqrt t\wedge R_0)/(2L).$
By Proposition \ref{P:3.3} and the parabolic Harnack inequality in Proposition \ref{P:3.1}, there exists $c_2=c_2(d, R_0, \Lambda_0, \phi, T)$ such that for $t\in (0, T)$ and $x\in D\cap B(z_0, (R_0\wedge\sqrt{t})/(8L)),$
\begin{equation}\label{e:3.21}
\begin{aligned}
\P_x(\tau_{U_{\sqrt t\wedge R_0}}>t)&=\int_{U_{\sqrt t\wedge R_0}} p_{U_{\sqrt t\wedge R_0}}(t, x, y)\,dy\\
&=\int_{U_{\sqrt t\wedge R_0}\cap \{y: \delta_{U_{\sqrt t\wedge R_0}}(y)\leq (\sqrt t\wedge R_0)/(16L)\}} p_{U_{\sqrt t\wedge R_0}}(t, x, y)\,dy\\
&\quad+\int_{U_{\sqrt t\wedge R_0}\cap \{y: \delta_{U_{\sqrt t\wedge R_0}}(y)> (\sqrt t\wedge R_0)/(16L)\}} p_{U_{\sqrt t\wedge R_0}}(t, x, y)\,dy\\
&\leq c_2 p_{U_{\sqrt t\wedge R_0}}(2t, x, y_0)(\sqrt t \wedge R_0)^{d}.
\end{aligned}
\end{equation}
Hence it follows from \eqref{e:3.11} and \eqref{e:3.21} that
\begin{equation}\label{e:3.26}
\P_x(\tau_D>t) \leq c_1c_2 p_{U_{\sqrt t\wedge R_0}}(2t, x, y_0)(\sqrt t \wedge R_0)^{d}, \quad x\in D\cap B(z_0, (R_0\wedge\sqrt{t})/(8L)).
\end{equation}
The desired conclusion is obtained.

\qed

In the following, we prove the doubling property of the survival probability for $X$ in $D$.

\begin{prp}\label{P:3.4}
Let $D$ be a Lipschitz open set with characteristics $(R_0, \Lambda_0).$
For each $T>0,$ there exists $C=C(d, \phi, R_0, \Lambda_0, T)>0$ such that for $t\in (0, T),$
\begin{equation}\label{e:3.24}
\P_x(\tau_D>t)\leq C\P_x(\tau_D>2t) \quad \hbox{for} \quad x\in D.
\end{equation}
\end{prp}

\pf Let $T>0$ and $t\in (0, T).$ If $x\in D$ with $\delta_D(x)> (\sqrt t\wedge R_0)/(8L),$ then by Lemma \ref{L:2.7}, there exists $c_1$ such that  for $t\in (0, T),$
$$\P_x(\tau_D>2t)\geq \P_x(\tau_{B(x, \frac{\sqrt t\wedge R_0}{16L})}>2t)\geq c_1\geq c_1\P_x(\tau_D>t),$$
Thus it suffices to prove \eqref{e:3.24}  for $x\in D$ with $\delta_D(x)\leq (\sqrt t\wedge R_0)/(8L).$

Let $x\in D$ with $\delta_D(x)\leq (\sqrt t\wedge R_0)/(8L).$
Let $z_x\in \partial D$ such that $|x-z_x|=\delta_D(x).$
For the simplicity of notation, we denote by $U_{\sqrt t\wedge R_0}:=U_{\sqrt t\wedge R_0}(z_x).$
By Proposition \ref{P:3.5}, there exists $c_2=c_2(d, \phi, R_0, \Lambda_0, T)>0$ such that
\begin{equation}\label{e:3.26}
\P_x(\tau_D>t) \leq c_2 p_{U_{\sqrt t\wedge R_0}}(2t, x, y_0)(\sqrt t \wedge R_0)^{d},
\end{equation}
where $y_0$ is a point of $D\cap \partial B(z, (\sqrt{t}\wedge R_0)/(2L))$ with $\delta_D(y_0)\geq \kappa (\sqrt{t}\wedge R_0)/(2L).$
On the other hand, by the parabolic Harnack inequality in Proposition \ref{P:3.1}, there exists $c_3=c_3(d, R_0, \Lambda_0, \phi, T)>0$ such that for $t\in (0, T),$
\begin{equation}\label{e:3.22}
\begin{aligned}
\P_x(\tau_D>4t)&\geq \P_x(\tau_{U_{\sqrt t\wedge R_0}}>4t)\\
&\geq \int_{B(y_0, \kappa (\sqrt t\wedge R_0)/(4L))} p_{U_{\sqrt t\wedge R_0}}(4t, x, y)\,dy\\
&\geq c_3p_{U_{\sqrt t\wedge R_0}}(2t, x, y_0) (\sqrt t \wedge R_0)^{d}.
\end{aligned}
\end{equation}
By combining \eqref{e:3.26} and \eqref{e:3.22},  for  $x\in D$ with $\delta_D(x)\leq (\sqrt t\wedge R_0)/(8L)$ and $t\in (0, T),$
$$\P_x(\tau_D>t)\leq \dfrac{c_2}{c_3}\P_x(\tau_D>4t)\leq \dfrac{c_2}{c_3}\P_x(\tau_D>2t).$$
The proof is complete.

\qed

\section{Proof of Theorem \ref{T2}}
\subsection{Upper bounds in Theorem \ref{T2}}

By an argument similar to that of \cite[Lemma 3.1]{CKS7} and the strong Markov property of $X,$ we have the following Lemma.

\begin{lem}\label{L:4.1}
Suppose that $U_1, U_3, E$ are open subsets of $\R^d$ with $U_1, U_3\subset E$ and ${\rm dist}(U_1, U_3)>0.$
Let $U_2:=E\setminus (U_1\cup U_3).$ If $x\in U_1$ and $y\in U_3,$ then for every $t>0,$
$$\begin{aligned}
p_E(t, x, y)
&\leq \E_x(p_E(t-\tau_{U_1}, X_{\tau_{U_1}}, y); \tau_{U_1}<t, X_{\tau_{U_1}}\in U_2)\\
&\quad+\E_x(p_E(t-\tau_{U_1}, X_{\tau_{U_1}}, y); \tau_{U_1}<t, X_{\tau_{U_1}}\in U_3)\\
&\leq \P_x(X_{\tau_{U_1}}\in U_2)\left(\sup_{s<t, z\in U_2}p_E(s, z, y)\right)\\
&+\int_0^t \P_x(\tau_{U_1}>s)\P_y(\tau_E>t-s)\,ds \left(\sup_{u\in U_1, z\in U_3} j(|u-z|)\right).
\end{aligned}$$
\end{lem}

\begin{thm}\label{T:5.1'}
Let $D$ be a Lipschitz open set in $\R^d$ with characteristics $(R_0, \Lambda_0).$

(i) For each $T>0,$ there exists $C_1=C_1(d, \phi, R_0, \Lambda_0, T)>0$  such that for any  $t\in (0, T)$ and $x, y\in D$ with $|x-y|\leq \sqrt t,$
\begin{equation}\label{e:4.3'}
p_D(t, x, y)\leq C_1t^{-d/2}\P_x(\tau_D>t)\P_y(\tau_D>t).
\end{equation}

(ii) There exist $C_2=C_2(d, \phi, R_0, \Lambda_0, T)>0$ and $C_3=C_3(d, \phi)>0$   such that for any  $t\in (0, T)$ and  $x, y\in D$ with $|x-y|> \sqrt t,$
\begin{equation}\label{e:4.4'}
p_D(t, x, y)\leq C_2\P_x(\tau_D>t)p(t, C_3x, C_3y)\P_y(\tau_D>t)+C_2\int_0^t \P_x(\tau_D>s)\P_y(\tau_D>t-s)\,ds \cdot j(|x-y|).
\end{equation}
\end{thm}

\pf Let $T>0$ and $t\in (0, T).$ We deal with two cases separately.
\smallskip

(i) Suppose $x, y\in D$ with $|x-y|\leq 8\sqrt{t}.$ By Theorem \ref{T1}, there exists $c_1=c_1(d, \phi, T)>0$ such that $p(t, w, z)\leq c_1t^{-d/2}$ for $t\in (0, T)$ and $w, z\in\R^d.$
By the semigroup property and the symmetry of $p_D(t, x, y)$ in $x$ and $y$,
\begin{equation}\label{e:4.1'}\begin{aligned}
\sup_{w\in D}p_D(t, x, w)&=\sup_{w\in D}\int_D p_D(t/2, x, z)p_D(t/2, z, w)\,dz\\
&\leq \sup_{z, w\in D} p(t/2, z, w)\int_D p_D(t/2, x, z)\,dz\\
&\leq c_1 \P_x(\tau_D>t/2) t^{-d/2}\\
\end{aligned}\end{equation}
Thus by \eqref{e:4.1'} and Proposition \ref{P:3.4}, for $t\in (0, T)$ and $|x-y|\leq 8t^{1/2},$
\begin{equation}\label{e:4.3n}\begin{aligned}
p_D(t, x, y)&=\int_D p_D(t/2, x, w)p_D(t/2, w, y)\,dw\\
&\leq \sup_{w\in D} p_D(t/2, x, w)\int_D p_D(t/2, w, y)\,dw\\
&\leq c_1 \P_x(\tau_D>t/4) (t/2)^{-d/2}\P_y(\tau_D>t/2)\\
&\leq c_2 t^{-d/2}\P_x(\tau_D>t) \P_y(\tau_D>t),
\end{aligned}\end{equation}
where $c_2=c_2(d, R_0, \Lambda_0, \phi, T)>0.$
Consequently, \eqref{e:4.3'} holds.

\smallskip

(ii) Suppose $x, y\in D$ with $|x-y|> 8\sqrt t.$
We first prove that there exist $c_3=c_3(d, \phi, R_0, \Lambda_0, T)$ and  $c_4=c_4(d, \phi)$ such that for any $t\in (0, T)$ and  $x, y\in D$ with
$|x-y|> 4\sqrt t,$
\begin{equation}\label{e:4.2'}
p_D(t, x, y)\leq c_3 \P_x(\tau_D>t)p(t, c_4x, c_4y)+ c_3\E_x (\tau_D \wedge t)\cdot j(|x-y|).
\end{equation}
For  $x\in D$ with $\delta_D(x)> (\sqrt t\wedge R_0)/(8L),$  by Lemma \ref{L:2.7}, there exists a constant $c_5=c_5(d, \phi, R_0, L, T)\in (0, 1)$ such that for $t\in (0, T),$
$$\P_x(\tau_D>t)\geq \P_x(\tau_{B(x, (\sqrt t\wedge R_0)/(8L))}>t)\geq c_5.$$
Note that $p_D(t, x, y)\leq p(t, x, y).$
Hence \eqref{e:4.2'} holds  in this case.

Next we shall  prove \eqref{e:4.2'}  for  $|x-y|> 4\sqrt t$ with $\delta_D(x)\leq (\sqrt t\wedge R_0)/(8L).$
Let $t_0:=\sqrt t\wedge R_0.$
Let $z_x\in \partial D$ such that $\delta_D(x)=|x-z_x|.$ Define
$$U_1:= U_{t_0}(z_x)$$
and
$$U_3:=\{z\in D: |z-x|\geq |x-y|/2\}, \quad U_2:=D\setminus (U_1\cup U_3).$$
Recall that $D\cap B(z_x, t_0/L)\subset U_1\subset D\cap B(z_x, t_0).$
Let $x_0$ be a point of $D\cap\partial B(z_x, t_0/(2L))$ with $\delta_D(x_0)\geq \kappa t_0/(2L).$
By  Proposition \ref{P:2.2} (i), there exists $c_6=c_6(d, \phi, R_0, \Lambda_0, T)$ such that
\begin{equation}\label{e:4.1}
\P_x(X_{\tau_{U_1}}\in D)\leq c_6\P_x(\tau_D>t).
\end{equation}
By Theorem \ref{T1}, there exist $c_7=c_7(d, \phi, T)>1$ and $c_8=c_8(d, \phi)>0$ such that for any $s\in (0, T)$ and $z\in\R^d,$
\begin{equation}\label{e:4.6'}
 p(s, z, y) \leq c_7\, s^{-d/2}\wedge \big(s^{-d/2}e^{-|z-y|^2/(c_8s)}+sj(|z-y|)\big).
\end{equation}
Observe that for any
$\delta>0$, the function $f(s):=s^{-d/2}e^{-\delta/s}$ is increasing
 on the interval $(0, 2\delta /d]$.
 Note that $|z-y|\geq |x-y|-|z-x|\geq |x-y|/2$ for $z\in U_2.$
By \eqref{e:4.6'} and the condition \eqref{e:1.1}, for $t  \leq |x-y|^2/(2dc_8)$,
\begin{equation}\label{e:4.7'}
\begin{aligned}
\sup_{s\le t\leq |x-y|^2/(2dc_8),\, z\in U_2} p(s, z, y)
& \leq
c_7 \sup_{s\le t\leq |x-y|^2/(2dc_8)} s^{-d/2}e^{-|x-y|^2/(4c_8s)} +
c_9 t
j(|x-y|)\\
&\leq
c_7t^{-d/2}e^{-|x-y|^2/(4c_8t)} + c_9 tj(|x-y|).
\end{aligned}\end{equation}
For $|x-y|^2/(2dc_8)\leq s<t\leq |x-y|^2/16,$ we have $8t/(dc_8)\leq s<t.$ Thus by \eqref{e:4.6'} and the condition \eqref{e:1.1},
\begin{equation}\label{e:4.8}
\sup_{|x-y|^2/(2dc_8)\leq s\leq t, \, z\in U_2} p(s, z, y)\leq c_{10}t^{-d/2}e^{-|x-y|^2/(4c_8t)} + c_9 tj(|x-y|).
\end{equation}
Hence it follows from \eqref{e:4.7'}, \eqref{e:4.8} and Theorem \ref{T1} that there exist $c_{11}=c_{11}(d, \phi, T)>1$ and $c_{12}=c_{12}(d, \phi)>0$ such that for $t\in (0, T)$ and $|x-y|> 4\sqrt t,$
\begin{equation}\label{e:4.2}
\sup_{0< s\leq t, \, z\in U_2} p(s, z, y)\leq c_{11}p(t, c_{12}x, c_{12}y).
\end{equation}
 Note that $|u-z|\geq |x-z|-|x-u|\geq |x-y|/4$ for $u\in U_1$ and $z\in U_3,$ by the condition \eqref{e:1.1},
\begin{equation}\label{e:4.3}
\sup_{u\in U_1, z\in U_3}j(|u-z|)=\sup_{u\in U_1, z\in U_3}\dfrac{1}{|u-z|^d \phi(|u-z|)}\leq c_{13}\dfrac{1}{|x-y|^d \phi(|x-y|)}= c_{13}j(|x-y|).
\end{equation}
By Lemma \ref{L:4.1}, \eqref{e:4.1}, \eqref{e:4.2} and \eqref{e:4.3}, for $|x-y|> 4\sqrt t$ with $\delta_D(x)\leq (\sqrt t\wedge R_0)/(8L),$
\begin{equation}\label{e:4.4}\begin{aligned}
p_D(t, x, y)
&\leq \P_x(X_{\tau_{U_1}}\in D)\left(\sup_{s<t, z\in U_2}p_D(s, z, y)\right)\\
&+\int_0^t \P_x(\tau_{U_1}>s)\P_y(\tau_D>t-s)\,ds \cdot\left(\sup_{u\in U_1, z\in U_3} j(|u-z|)\right)\\
&\leq c_6c_{11}\P_x(\tau_D>t) p(t, c_{12}x, c_{12}y)+c_{13}\int_0^t \P_x(\tau_D>s)\,ds \cdot j(|x-y|)\\
&\leq c_6c_{11}\P_x(\tau_D>t) p(t, c_{12}x, c_{12}y)+c_{13}\E_x (\tau_D\wedge t) \cdot j(|x-y|).
\end{aligned}\end{equation}
Hence \eqref{e:4.2'} holds.

\smallskip

In the following, we use \eqref{e:4.2'} to prove that \eqref{e:4.4'} holds.
Suppose $|x-y|> 8\sqrt t.$
Let $D_y:=\{z\in D: |z-y|\leq |x-y|/2\}.$ We have
\begin{equation}\label{e:4.8'}\begin{aligned}
p_D(t, x, y)&=\int_{D\setminus D_y} p_D(t/2, x, z)p_D(t/2, z, y)\,dz+\int_{D_y} p_D(t/2, x, z)p_D(t/2, z, y)\,dz\\
&\leq \P_x(\tau_D>t/2)\sup_{z\in D\setminus D_y}p_D(t/2, z, y)+\P_y(\tau_D>t/2)\sup_{z\in D_y}p_D(t/2, x, z).
\end{aligned}\end{equation}
Since $|z-y|\geq |x-y|/2\geq 4\sqrt t$ for $z\in D\setminus D_y,$ then by \eqref{e:4.2'}, Proposition \ref{P:3.4} and the condition \eqref{e:1.1},  for $z\in D\setminus D_y,$
\begin{equation}\label{e:4.9'}\begin{aligned}
p_D(t/2, z, y)&\leq c_3\P_y(\tau_D>t/2) p(t, c_4z, c_4y)+c_3\E_y (\tau_D\wedge t/2) \cdot j(|z-y|)\\
&\leq c_{14}\P_y(\tau_D>t) p(t,c_{15} x, c_{15}y)+c_{14}\E_y (\tau_D\wedge t) \cdot j(|x-y|).
\end{aligned}\end{equation}
Similarly, since $|x-z|\geq |x-y|-|y-z|\geq |x-y|/2\geq 4\sqrt t$ for $z\in D_y,$  we have
\begin{equation}\label{e:4.10}
\sup_{z\in D_y}p_D(t/2, x, z)\leq c_{14}\P_x(\tau_D>t)p(t, c_{15}x, c_{15}y)+c_{14} \E_x(\tau_D\wedge t)j(|x-y|).
\end{equation}
By \eqref{e:4.8'}-\eqref{e:4.10} and Proposition \ref{P:3.4}, for $|x-y|> 8\sqrt t,$
\begin{equation}\label{e:4.12}\begin{aligned}
p_D(t, x, y)&\leq c_{14}\P_x(\tau_D>t/2)\cdot \Big[ p(t, c_{15}x, c_{15}y)\P_y(\tau_D>t)+\E_y(\tau_D\wedge t)j(|x-y|)\Big]\\
&\quad+c_{14}\P_y(\tau_D>t/2)\cdot \Big[ p(t, c_{15}x, c_{15}y)\P_x(\tau_D>t)+\E_x(\tau_D\wedge t)j(|x-y|)\Big]\\
&\leq c_{16}\P_x(\tau_D>t) p(t, c_{15}x,c_{15}y)\P_y(\tau_D>t)\\
&\quad +c_{16}\Big[\P_x(\tau_D>t)\E_y (\tau_{D}\wedge t)+\P_y(\tau_D>t)\E_x (\tau_{D}\wedge t)\Big]\cdot j(|x-y|).
\end{aligned}\end{equation}
We denote $f\overset{c}{\asymp} g$ if $c^{-1}f\leq g\leq c f$ in the common domain of $f$ and $g.$
By Proposition \ref{P:3.4}, there exists $c_{17}=c_{17}(d, \phi, R_0, \Lambda_0, T)>1$ such that for  $t\in (0, T)$ and $x, y\in D,$
\begin{equation}\label{e:4.13}\begin{aligned}
&\int_0^t \P_x(\tau_D>s)\P_y(\tau_D>t-s)\,ds \\
&=\int_0^{t/2} \P_x(\tau_D>s)\P_y(\tau_D>t-s)\,ds +\int_{t/2}^t \P_x(\tau_D>s)\P_y(\tau_D>t-s)\,ds\\
&\overset{c_{17}}{\asymp} \P_y(\tau_D>t)\int_0^{t/2} \P_x(\tau_D>s)\,ds+\P_x(\tau_D>t)\int_0^{t/2}\P_y(\tau_D>s)\,ds\\
&\overset{c_{17}}{\asymp} \P_y(\tau_D>t)\E_x (\tau_D\wedge \frac{t}{2})+\P_x(\tau_D>t) \E_y (\tau_D\wedge \frac{t}{2}).
\end{aligned}\end{equation}
Hence, for  $t\in (0, T)$ and $x, y\in D,$
\begin{equation}\label{e:4.13n}
\begin{aligned}
&(2c_{17})^{-1}\Big[ \P_y(\tau_D>t)\E_x (\tau_D\wedge t)+\P_x(\tau_D>t) \E_y (\tau_D\wedge t)\Big]\\
&\leq
\int_0^t \P_x(\tau_D>s)\P_y(\tau_D>t-s)\,ds\\
&\leq
c_{17} \Big[\P_y(\tau_D>t)\E_x (\tau_D\wedge t)+\P_x(\tau_D>t) \E_y (\tau_D\wedge t)\Big].
\end{aligned}
\end{equation}
Consequently, \eqref{e:4.4'}  holds for $x, y\in D$ with $|x-y|> 8\sqrt t$ by \eqref{e:4.12} and \eqref{e:4.13n}.
By Theorem \ref{T1}, there exists $c_{18}=c_{18}(d, \phi, T)>0$ such that $p(t, x, y)\geq c_{18}t^{-d/2}$ for $t\in (0, T)$ and $|x-y|\leq 8t^{1/2}.$
Then by \eqref{e:4.3n}, it is easy to see that \eqref{e:4.4'}  holds for $x, y\in D$ with $\sqrt t<|x-y|\leq 8\sqrt t.$
The proof is complete.
\qed

\subsection{Lower bounds in Theorem \ref{T2}}

\begin{lem}\label{L:4.4}
For each $T>0$,  there exists $C=C(d, \phi, R_0, \Lambda_0, T)>0$ such that for any  $t\in (0, T)$ and
$x, y\in D$ with $|x-y|\geq \sqrt t\wedge R_0,$
\begin{equation}\label{e:4.16}
p_D(t, x, y) \geq C\int_0^t \P_x(\tau_D>s)\P_y(\tau_D>t-s)\,ds \cdot j(|x-y|).
\end{equation}
\end{lem}

\pf For $t\in (0, T),$ let $t_0:=(\sqrt t\wedge R_0)/4.$ Let $x, y\in D$ with $|x-y|\geq 4 t_0.$

(i) We first consider the case for $x, y\in D$ with $|x-y|>4t_0$ and $\delta_D(x)\leq t_0/(8L), \delta_D(y)\leq t_0/(8L).$
Let $z_x$ and $z_y$  be  points  of $\partial D$ such that $|x-z_x|=\delta_D(x)$ and $|y-z_y|=\delta_D(y).$
Let $U_1:=U_{t_0}(z_x)$ and $U_3= U_{t_0}(z_y).$

By the strong Markov property of $X$ and the L\'evy system formula \eqref{e:2.2}, we have
\begin{equation}\label{e:4.26}\begin{aligned}
p_D(t, x, y)
&\geq \E_x[p_D(t-\tau_{U_1}, X_{\tau_{U_1}}, y); X_{\tau_{U_1}}\in U_3, \tau_{U_1}<t]\\
&= \int_0^{t} \left(\int_{U_1} p_{U_1}(s, x, w)\left(
\int_{U_3} j(|w-z|) p_D (t-s, z, y) dz\right) dw \right) ds \\
&\geq \int_0^t \P_x(\tau_{U_1}>s)\P_y(\tau_{U_3}>t-s)\,ds \cdot \left(\inf_{w\in U_1, z\in U_3} j(|w-z|)\right).
\end{aligned}\end{equation}
Note that
if
$w\in U_1$ and
$z\in U_3$, then
$$
|w-z| \le |x-w| +|x-y|+|y-z| \le |x-y| + 4t_0 \le 2|x-y|.
$$
By \eqref{e:4.26} and the condition  \eqref{e:1.1}, we have
\begin{equation}\label{e:4.17n}
p_D(t, x, y)
\geq c_1\int_0^t\P_x(\tau_{U_1}>s)\P_y(\tau_{U_3}>t-s)\,ds \cdot j(|x-y|).
\end{equation}
By Proposition \ref{P:2.2} (ii), there exists $c_2=c_2(d, \phi, R_0, \Lambda_0, T)$ such that for $t\in (0, T)$ and $s\in (0, t),$
\begin{equation}\label{e:4.28'}
\P_x(\tau_{U_1}>s)\geq c_2\P_x(\tau_D>s) \quad \hbox{and} \quad \P_y(\tau_{U_3}>s)\geq c_2\P_y(\tau_D>s).
\end{equation}
By  \eqref{e:4.17n} and \eqref{e:4.28'},  \eqref{e:4.16} holds for  $|x-y|>4t_0$ with $\delta_D(x)\leq t_0/(8L)$ and $\delta_D(y)\leq t_0/(8L).$

\smallskip

(ii) For $x, y\in D$ with $|x-y|>4t_0, \delta_D(x)> t_0/(8L)$ and $\delta_D(y)\leq t_0/(8L),$
let $B_x=B(x, t_0/(8L)).$
Then by an argument similar to that of \eqref{e:4.26} and \eqref{e:4.17n} but with $B_x$ in place of $U_1$,
\begin{equation}\label{e:4.21n}
p_D(t, x, y)\geq c_1\int_0^t \P_x(\tau_{B_x}>s)\P_y(\tau_{U_3}>t-s)\,ds \cdot j(|x-y|).
\end{equation}
By Lemma \ref{L:2.7}, there exists $c_3=c_3(d, \phi, T)$ such that for $t\in (0, T)$ and $s\in (0, t),$
\begin{equation}\label{e:4.27'}
\P_x(\tau_{B_x}>s)\geq \P_x(\tau_{B_x}>t)\geq c_3\geq c_3\P_x(\tau_D>s).
\end{equation}
Hence by \eqref{e:4.21n}, \eqref{e:4.27'} and the second inequality of \eqref{e:4.28'},  \eqref{e:4.16} holds for $|x-y|>4t_0$ with $\delta_D(x)\geq t_0/(8L)$ and $\delta_D(y)\leq t_0/(8L),$
By the symmetry of $p_D(t, x, y)$ in $x$ and $y,$ then  \eqref{e:4.16} holds for  $|x-y|>4t_0$  with $\delta_D(x)\leq t_0/(8L)$ and $\delta_D(y)> t_0/(8L).$

\smallskip

(iii) For $x, y\in D$ with $|x-y|>4t_0, \delta_D(x)> t_0/(8L)$ and $\delta_D(y)> t_0/(8L),$
let $B_x=B(x, t_0/(8L))$ and  $B_y=B(y, t_0/(8L)).$
By an argument similar to that of \eqref{e:4.26} and \eqref{e:4.17n} but  with $B_x$ and $B_y$ in place of $U_1$ and $U_3$,
\begin{equation}\label{e:4.23n}
p_D(t, x, y)\geq \int_0^t \P_x(\tau_{B_x}>s)\P_y(\tau_{B_y}>t-s)\,ds \cdot j(|x-y|).
\end{equation}
Then by  \eqref{e:4.27'},  \eqref{e:4.16} holds in this case.

The proof is complete.
\qed

\bigskip

The next lemma can be obtained by an argument  similar to that of \cite[Lemma 2.5]{CKS7}.
 It implies
that if $x$ and $y$ are in different components of $D$,
the jumping kernel component of the heat kernel
dominates the Gaussian component.

\begin{lem}\label{L:4.3'}
For any given positive constants $c_1, c_2, b$ and $T,$ there is a positive constant $C=C(d, c_1, c_2, b, T, \phi)$ such that
$$t^{-d/2}e^{-r^2/c_1t}\leq C(t^{-d/2}\wedge (tj(c_2r)))\qquad \hbox{for every } r\geq b \hbox{ and } t\in (0, T].$$
\end{lem}

\begin{thm}\label{T:5.2'}
Suppose that $D$ is a Lipschitz open set in $\R^d$ with characteristics $(R_0, \Lambda_0).$ Assume the path distance in each connected component of $D$ is comparable to the Euclidean distance with characteristic $\chi_1.$

(i) For each $T>0,$ there exists $C_1=C_1(d, \phi, R_0, \Lambda_0, \chi_1, T)>0$ such that for  any $t\in (0, T)$ and $x, y\in D$ with $|x-y|\leq \sqrt t,$
\begin{equation}\label{e:4.17'}
p_D(t, x, y)\geq C_1t^{-d/2}\P_x(\tau_D>t)\P_y(\tau_D>t).
\end{equation}

(ii)  There exist positive constants $C_2=C_2(d, \phi, R_0, \Lambda_0, \chi_1, T)$ and $C_3=C_3(d, \phi, \chi_1)$  such that for  any $t\in (0, T)$ and $x, y\in D$ with $|x-y|> \sqrt t,$
\begin{equation}\label{e:4.18}
p_D(t, x, y)\geq C_2\P_x(\tau_D>t)p(t, C_3x, C_3y)\P_y(\tau_D>t)+C_2\int_0^t \P_x(\tau_D>s)\P_y(\tau_D>t-s)\,ds\cdot j(|x-y|).
\end{equation}
\end{thm}

\pf Let $D$ be a Lipschitz open set in $\R^d$ with characteristics $(R_0, \Lambda_0).$
For each $t\in (0, T),$ let $t_0:=\sqrt t\wedge R_0.$
For $x, y\in D,$ let $z_x$ and $z_y$ be the points  of $\partial D$ such that $|x-z_x|=\delta_D(x)$ and $|y-z_y|=\delta_D(y).$
If $x\in D$ with $\delta_D(x)<t_0/(32L),$  let $A^{t_0}_x$ be a point of $D\cap\partial B(z_x, t_0/(2L))$ with $\delta_D(A^{t_0}_x)\geq \kappa t_0/(2L).$
Otherwise if $x\in D$ with $\delta_D(x)>t_0/(32L),$ we let $A^{t_0}_x=x.$
Similarly we can define $A^{t_0}_y$ for $y\in D.$

 By the semigroup property of $p_D(t, x, y)$, for $ t>0$ and $x, y\in D,$
\begin{equation}\label{e:4.17}
p_D(t, x, y)\geq \int_{B(A^{t_0}_y, \kappa t_0/(64L))}\int_{B(A^{t_0}_x, \kappa t_0/(64L))}p_D(t/3, x, u)p_D(t/3, u, v)p_D(t/3, v, y)\,du\,dv.
\end{equation}
We first prove that there exists $c_1=c_1(d, \phi, R_0, \Lambda_0, T)$ such that for $t\in (0, T),$
\begin{equation}\label{e:4.6}
\int_{B(A^{t_0}_x, \kappa t_0/(64L))}p_D(t/3, x, u)\,du\geq  c_1\P_x(\tau_D>t), \quad x\in D.
\end{equation}
In fact, if $x\in D$ with $\delta_D(x)>t_0/(32L),$ note that $t_0^2\leq t\leq (1\vee \frac{T}{R_0^2})t_0^2$ and $A^{t_0}_x=x,$   there exists $c_2=c_2(d, \phi, R_0, \Lambda_0, T)\in (0, 1)$ such that
\begin{equation}\label{e:4.4}
\begin{aligned}
\int_{B(A^{t_0}_x, \kappa t_0/(64L))}p_D(t/3, x, u)\,du&\geq
 \int_{B(x, \kappa t_0/(64L))}p_{B(x, \kappa t_0/(64L))}(t/3, x, u)\,du\\
 &= \P_x(\tau_{B(x, \kappa t_0/(64L))}>t/3)\geq c_2\geq c_2 \P_x(\tau_D>t),
 \end{aligned}
\end{equation}
where the second inequality is due to  Lemma \ref{L:2.7}.
For $x\in D$ with $\delta_D(x)\leq t_0/(32L),$  by the parabolic Harnack inequality in Proposition \ref{P:2.1}, there exists $c_3=c_3(d, \phi, R_0, \Lambda_0, T)>0$ such that
\begin{equation}\label{e:4.24'}
p_D(t/3, x, u)\geq c_3 p_D(t/6, x, A^{t_0}_x) \quad \hbox{for } u\in B(A^{t_0}_x, \kappa t_0/(64L)).
\end{equation}
On the other hand, it follows from Proposition \ref{P:3.5} that there exists $c_4=c_4(d, \phi, R_0, \Lambda_0, T)$ such that
\begin{equation}\label{e:4.23'}
\P_x(\tau_D>t/12) \leq c_4 p_D(t/6, x, A^{t_0}_x)t_0^d.
\end{equation}
Then by combining \eqref{e:4.24'} and \eqref{e:4.23'},
\begin{equation}\label{e:4.5'}
\begin{aligned}
\int_{B(A^{t_0}_x, \kappa t_0/(64L))}p_D(t/3, x, u)\,du
&\geq c_5 p_D(t/6, x, A^{t_0}_x) t_0^d\\
&\geq c_5c_4^{-1}\P_x(\tau_D>t/12)\\
&\geq c_5c_4^{-1}\P_x(\tau_D>t).
\end{aligned}
\end{equation}
Hence \eqref{e:4.6} is obtained by  \eqref{e:4.4} and \eqref{e:4.5'}.

\smallskip

(a) Suppose $x, y\in D$ with $|x-y|<t_0.$ Since $D$ is a Lipschitz open set with characteristics $(R_0, \Lambda_0),$  the distance between two distinct connected components of D is at least $R_0.$ Note that $t_0\leq R_0.$ In this case $x$ and $y$ are in the same connected component of $D.$
Hence $B(A^{t_0}_x, \kappa t_0/(64L))$ and $B(A^{t_0}_y, \kappa t_0/(64L))$ are contained in a connected component of $D.$
Note that $|u-v|\leq |u-A^{t_0}_x|+|A^{t_0}_x-A^{t_0}_y|+|A^{t_0}_y-v|\leq c_6t_0$ for $(u, v)\in B(A^{t_0}_x, \kappa t_0/(64L))\times B(A^{t_0}_y, \kappa t_0/(64L)),$
there is a rectifiable curve $l$ in $U$ which connects $u$ to $v$ such that the length of $l$ is less than or equal to $ c_6\chi_1 t_0.$
Thus by \eqref{e:2.9} and the standard  chain argument (see e.g. \cite[Lemma 2.1]{Zh}), there exists $c_7=c_7(d, \phi, \chi_1, T)>0$ such that
$$p_D(t, u, v)\geq c_7t^{-d/2} \quad \mbox{for} \quad (u, v)\in B(A^{t_0}_x, \kappa t_0/(64L))\times B(A^{t_0}_y, \kappa t_0/(64L)). $$
By this inequality and \eqref{e:4.6}, for $x, y\in D$ with $|x-y|<t_0,$
\begin{equation}\label{e:4.32'}\begin{aligned}
&p_D(t, x, y)\\
\geq &\int_{B(A^{t_0}_y, \kappa t_0/(64L))}\int_{B(A^{t_0}_x, \kappa t_0/(64L))}p_D(t/3, x, u)p_D(t/3, u, v)p_D(t/3, v, y)\,du\,dv\\
\geq &c_8\P_x(\tau_D>t)\P_y(\tau_D>t)t^{-d/2},
\end{aligned}\end{equation}
where $c_8:=c_1^2c_7.$ Hence \eqref{e:4.17'} holds for $|x-y|\leq \sqrt t\wedge R_0.$

(b) Suppose $x, y\in D$ with $|x-y|\geq t_0.$ By \eqref{e:4.6},
\begin{equation}\label{e:4.20'}\begin{aligned}
&p_D(t, x, y)\\
\geq &\int_{B(A^{t_0}_y, \kappa t_0/(64L))}\int_{B(A^{t_0}_x, \kappa t_0/(64L))}p_D(t/3, x, u)p_D(t/3, u, v)p_D(t/3, v, y)\,du\,dv\\
\geq &c_1^2\P_x(\tau_D>t)\P_y(\tau_D>t)
\inf_{(u, v)\in B(A^{t_0}_x, \kappa t_0/(64L))\times B(A^{t_0}_y, \kappa t_0/(64L))}p_D(t/3, u, v).
\end{aligned}\end{equation}
If $x$ and $y$ are in the same connected component of $D,$ then $B(A^{t_0}_x, \kappa t_0/(64L))$ and $B(A^{t_0}_y, \kappa t_0/(64L))$ are contained in a connected component of $D.$
Note that  $|u-v|\leq |u-A^{t_0}_x|+|A^{t_0}_x-A^{t_0}_y|+|A^{t_0}_y-v|\leq c_9|x-y|$ for $(u, v)\in B(A^{t_0}_x, \kappa t_0/(64L))\times B(A^{t_0}_y, \kappa t_0/(64L)).$
Hence, by \eqref{e:2.9} and the standard chain argument (see e.g. \cite[Lemma 2.1]{Zh}),  there exist $c_k=c_k(d, \phi, \chi_1)>0, k=10, 11$ such that
\begin{equation}\label{e:4.21'}
\begin{aligned}
\inf_{(u, v)\in B(A^{t_0}_x, \kappa t_0/(64L))\times B(A^{t_0}_y,\kappa  t_0/(64L))}p_D(t/3, u, v)
\geq c_{10} t^{-d/2} \exp(-|x-y|^2/(c_{11}t)).
\end{aligned}
\end{equation}
Hence by \eqref{e:4.20'} and \eqref{e:4.21'}, if $x$ and $y$ are in the same connected component of $D,$ then
\begin{equation}\label{e:4.34'}
p_D(t, x, y)\geq c_1^2c_{10}\P_x(\tau_D>t)\P_y(\tau_D>t)t^{-d/2} \exp(-|x-y|^2/(c_{11}t)).
\end{equation}

On the other hand, by Lemma \ref{L:4.4}, for each $T>0$,  there exists $c_{12}=c_{12}(d, \phi, R_0, \Lambda_0, T)>0$ such that for any  $t\in (0, T)$ and any
$x, y\in D$ with $|x-y|\geq t_0,$
\begin{equation}\label{e:4.34}
\begin{aligned}
p_D(t, x, y) &\geq c_{12}\int_0^t \P_x(\tau_D>s)\P_y(\tau_D>t-s)\,ds \cdot j(|x-y|)\\
&\geq c_{12} \P_x(\tau_D>t)\P_y(\tau_D>t) \cdot t j(|x-y|).
\end{aligned}
\end{equation}
Note that  the distance between two distinct connected components of D is at least $R_0.$
By Lemma \ref{L:4.3'}, there exists $c_{13}>0$ such that $t^{-d/2} \exp(-|x-y|^2/c_{11}t)\leq c_{13} tj(|x-y|)$ for $x$ and $y$ in different components of $D.$
Hence by combining  \eqref{e:4.34'} and  \eqref{e:4.34}, for any $x, y\in D$ with $|x-y|\geq t_0,$
\begin{equation}\label{e:4.7}
\begin{aligned}
p_D(t, x, y)&\geq c_{14}\P_x(\tau_D>t)\P_y(\tau_D>t)\cdot \big(t^{-d/2} \exp(-|x-y|^2/c_{11}t)+tj(|x-y|)\big)\\
&\geq  c_{15}\P_x(\tau_D>t)\P_y(\tau_D>t)p(t, c_{16}x, c_{16}y),
\end{aligned}
\end{equation}
where the last inequality is due to Theorem \ref{T1}.
Finally,  \eqref{e:4.7} together with Lemma \ref{L:4.4} yields \eqref{e:4.18} holds for $|x-y|>\sqrt t\wedge R_0$.
Hence \eqref{e:4.18} holds for $|x-y|>\sqrt t$.

Note that by Theorem \ref{T1}, there exists $c_{17}$ such that $p(t, c_{16}x, c_{16}y)\geq c_{17}t^{-d/2}$ for $|x-y|<\sqrt t.$
Combining this with  \eqref{e:4.7} and \eqref{e:4.32'}, \eqref{e:4.17'} holds for $|x-y|\leq \sqrt t.$ The proof is complete.

\qed

\medskip

{\bf Proof of Theorem \ref{T2}}
 Theorem \ref{T2} (i) and (ii) are obtained by Theorems \ref{T:5.1'} and \ref{T:5.2'}.
In the following, we prove Theorem \ref{T2} (iii).

Suppose $D$ is a bounded Lipschitz open set with characteristics $(R_0, \Lambda_0).$  Note that $|x-y|\leq {\rm diam}(D)$ for $x, y\in D.$
By virtue of Theorem \ref{T1}, there exists $c_1=c_1(d, \phi, {\rm diam}(D))>1$ such that $c_1^{-1}\leq p(1, x, y)\leq c_1 $ for $x, y\in D.$
Then by an argument similar to that of  \eqref{e:4.3n} with $t=1$, there exists $c_2=c_2(d, \phi, R_0, \Lambda_0, {\rm diam}(D))$ such that
$$p_D(1, x, y)\leq c_2\P_x(\tau_D>1)\P_y(\tau_D>1), \quad x, y\in D.$$
By Theorem \ref{T2}, it is easy to see that there exists $c_3=c_3(d, \phi, R_0, \Lambda_0, {\rm diam}(D))$ such that
$$p_D(1, x, y)\geq c_3\P_x(\tau_D>1)\P_y(\tau_D>1), \quad x, y\in D.$$
 By applying the two inequalities above and an argument  similar to that  of \cite[Theorem 1.3 (iii)]{CKS6}, the desired result is obtained.

\qed

\section{Proof of Theorem  \ref{T3} and Corollary \ref{C2}}

In this section, we shall prove Theorem  \ref{T3} and Corollary \ref{C2}.
Proposition \ref{L:4.5} below  is the key step in the proof of Theorem \ref{T3}.
In the following, for the simplicity of notation, for each open set $D, z\in \partial D$ and $r>0,$ we denote $D_r(z):=D\cap B(z, r).$

\begin{prp}\label{L:4.5}
Suppose that $D$ is a Lipschitz open set in $\R^d$ with characteristics $(R_0, \Lambda_0).$
 There exists $C=C(d, \phi, R_0, \Lambda_0, T)>0$ such that for  any $z_0\in \partial D, t\in (0, T)$ and $x\in D_{(\sqrt t\wedge R_0)/4}(z_0),$
 \begin{equation}\label{e:4.19}
 C^{-1}G_{D_{\sqrt t\wedge R_0}(z_0)}(x, x_t)t^{(d-2)/2}\leq \P_x(\tau_D>t)\leq CG_{D_{\sqrt t\wedge R_0}(z_0)}(x, x_t)t^{(d-2)/2},
 \end{equation}
 where $x_t$ is a point of $D\cap \partial B(z_0, (\sqrt t\wedge R_0)/2)$ with $\delta_D(x_t)\geq \kappa(\sqrt t\wedge R_0)/2.$
\end{prp}

\pf Let $t\in (0, T)$ and define $t_0:=\sqrt t\wedge R_0.$
Then $t_0^2\leq t\leq (1\vee \frac{T}{R_0^2})t_0^2.$
  For $z_0\in \partial D$ and $x\in D_{t_0/4}(z_0),$
 let  $x_t$ be a point of $D\cap \partial B(z_0, t_0/2)$ with $\delta_D(x_t)\geq \kappa t_0/2.$
 By an argument similar to that of Lemma \ref{L:2.6} with $D_{t_0}(z_0)$ in place of $U_{t_0}(z_0),$ there exists $c_1=c_1(d, \phi, \kappa)>0$ such that
$$
\P_x \left( X_{\tau_{D_{t_0}(z_0)\setminus \overline{B(x_t, \kappa t_0/4)}}}\in \overline{B(x_t, \kappa t_0/4)} \right)
\geq c_1 t_0^{d-2} G_{D_{t_0}(z_0)} (x, x_t) \quad  \hbox{for } x\in D_{t_0/4}(z).
$$
Let $F:=D_{t_0}(z)\setminus \overline{B(x_t, \kappa t_0/4)}.$
By  an argument similar to that of \eqref{e:2.10},
there exists $c_2=c_2(d, \phi, R_0, \Lambda_0, T)>0$ such that
 \begin{equation}\label{e:5.2}\begin{aligned}
\P_x(\tau_D>t)&\geq \P_x(\tau_{D_{t_0}(z_0)}>t)\\
&\geq \E_x[X_{\tau_F}\in B(x_t, \kappa t_0/4);\P_{X_{\tau_F}} (\tau_{B(X_{\tau_F}, \kappa t_0/4)}>t)]\\
&\geq \P_x[X_{\tau_F}\in B(x_t, \kappa t_0/4)] \inf_{y\in B(x_t, \kappa t_0/4)}\P_y(\tau_{B(y, \kappa t_0/4)}>t)\\
&\geq c_2G_{D_{t_0}(z_0)}(x, x_t)t_0^{d-2}\geq c_2(1\wedge (R_0/\sqrt T))^{d-2}G_{D_{t_0}(z_0)}(x, x_t)t^{(d-2)/2}.
\end{aligned}\end{equation}
That is, the lower bound of \eqref{e:4.19} holds.

In the following, we shall prove the upper bound of \eqref{e:4.19}.
Let $z_0\in \partial D$ and
$CS_{z_0}$ be an orthonormal coordinate system with its origin at $z_0$ such that
$ B(z_0, R_0)\cap D=\{y=(\wt {y}, y^d)\in B(0,R_0) \hbox{ in } CS_{z_0}: y_d>\Gamma_{z_0}(\wt {y})\}.$
For $x\in D_{t_0/4}(z_0),$ define $\rho_{z_0}(x):=x^d-\phi(\tilde x)$ in the coordinate $CS_{z_0}.$
Note that $G_{D_{t_0}(z_0)}(\cdot, x_t)$ is harmonic in $D_{t_0}(z_0) \setminus \{x_t\}$.
 For $x\in D_{t_0/4}(z_0)$ with $\rho_{z_0}(x)>t_0/(64L),$ by the standard chain argument, the elliptic Harnack inequality \eqref{e:2.5} and Lemma \ref{L:2.4'}, there exist $c_k=c_k(d, \phi, R_0, \Lambda_0, T), k=3, 4$ such that
$G_{D_{\sqrt t}(z_0)}(x, x_t)\geq c_3\inf_{y\in B(x_t, \kappa t_0/4)\setminus B(x_t, \kappa t_0/8)}G_{D_{\sqrt t}(z_0)}(y, x_t)\geq c_3\inf_{y\in B(x_t, \kappa t_0/4)}
G_{B(x_t, \kappa t_0/2)}(y, x_t)\geq c_4 t^{(2-d)/2}.$
Hence, for $x\in D_{t_0/4}(z_0)$ with $\rho_{z_0}(x)>t_0/(64L),$
$$G_{D_{t_0}(z_0)}(x, x_t)t^{(d-2)/2}\geq c_4\geq  c_4\P_x(\tau_D>t).$$
Thus it suffices to prove the upper bound of \eqref{e:4.19} for  $x\in D_{t_0/4}(z_0)$ with $\rho_{z_0}(x)\leq t_0/(64L).$
For  $x=(\tilde x, x^d)\in D_{t_0/4}(z_0)$ with $\rho_{z_0}(x)\leq t_0/(64L),$ let $z_x=(\tilde x, \Gamma_{z_0}(\tilde x)).$
Then $z_x\in \partial D\cap B(z_0, t_0/4).$
 By Propositions \ref{P:3.5} and \ref{P:3.1}, there exists
$c_5=c_5(d, \phi, R_0, \Lambda_0, T)>0$ such that
\begin{equation}\label{e:4.20}
\P_x(\tau_D>t)\leq c_5 t_0^{d} p_D(4t, x, x'_t),
\end{equation}
where $x'_t$ is a point of $D\cap \partial B(z_x, t_0/2)$ with $\delta_D(x'_t)\geq \kappa t_0/2.$

Let
$U_1:= U_{t_0/8}(z_x)$
and
$$U_3:=\{y\in D: |y-x|\geq |x-x'_t|/2\}, \quad U_2:=D\setminus (U_1\cup U_3).$$
Note that $D_{t_0/(8L)}(z_x)\subset U_1 \subset D_{t_0/8}(z_x).$
By Lemma \ref{L:4.1},
$$
 p_D(4t, x, x'_t)
 \leq \P_x(X_{\tau_{U_1}}\in D) \sup_{s<4t, u\in U_2}p(s, u, x'_t)+ \int_0^{4t} \P_x(\tau_{U_1}>s)\,ds\cdot \sup_{u\in U_1, y\in U_3}j(|u-y|).
 $$
 Note that  $|u-y|\geq |y-z_x|-|z_x-u|\geq |x-x'_t|/2-t_0/4\geq |x-x'_t|/32$ for $u\in U_1$ and $y\in U_3.$ By the condition \eqref{e:1.1} and an argument similar to that of \eqref{e:4.2} and \eqref{e:4.3},
 there exist $c_k=c_k(d, \phi, R_0, \Lambda_0, T), k=6, 7$ such that
 \begin{equation}\label{e:4.21}\begin{aligned}
 p_D(4t, x, x'_t)
 &\leq c_6\P_x(X_{\tau_{U_1}}\in D) p(4t, c_7x, c_7x'_t)+ c_6 \int_0^{4t} \P_x(\tau_{U_1}>s)\,ds\cdot j(|x-x'_t|)\\
 &\leq c_6\P_x(X_{\tau_{U_1}}\in D) p(4t, c_7x, c_7x'_t)+ c_6 \E_x\tau_{U_1} j(|x-x'_t|).
 \end{aligned}\end{equation}
By the L\'evy system formula \eqref{e:2.2} and the condition \eqref{e:1.1},
\begin{equation}\label{e:4.22'}\begin{aligned}
\P_x(X_{\tau_{U_1}}\in D)&\geq \P_x(X_{\tau_{U_1}}\in B(x'_t, \kappa t_0/16))\\
&=\E_x\int_0^{\tau_{U_1}}\int_{B(x'_t, \kappa t_0/16)} j(|X_s-y|)\,dy\,ds\\
&\geq c_8\E_x \tau_{U_1} t_0^d \cdot \inf_{u\in U_1, y\in B(x'_t, \kappa t_0/16)} j(|u-y|)\\
&\geq c_9t_0^{d}\E_x\tau_{U_1} j(|x-x'_t|).
\end{aligned}\end{equation}
By Theorem \ref{T1}, there exists $c_{10}=c_{10}(d, \phi, T)>0$ such that $p(4t, c_6x, c_6x'_t)\leq c_{10} t^{-d/2}\leq c_{10}t_0^{-d}$ for any $t\in (0, T).$
Hence by \eqref{e:4.21} and \eqref{e:4.22'},
\begin{equation}\label{e:4.22}
p_D(4t, x, x'_t)\leq c_6(c_{10}+c_9^{-1})t_0^{-d}\P_x(X_{\tau_{U_1}}\in D).
\end{equation}
Thus by combining \eqref{e:4.20} and \eqref{e:4.22},
\begin{equation}\label{e:5.6n}\begin{aligned}
\P_x(\tau_D>t)\leq c_{11}\P_x(X_{\tau_{U_1}}\in D),
\end{aligned}\end{equation}
where $c_{11}:=c_5c_6(c_{10}+c_9^{-1}).$
Let  $y_0$ be a point of $D\cap \partial B(z_x, t_0/(16L))$ with $\delta_D(y_0)\geq \kappa t_0/(16L).$
Note that $z_x\in \partial D\cap B(z_0, t_0/4)$, thus $U_1\subset D_{t_0}(z_0).$
By \eqref{e:5.6n} and Proposition \ref{P:2.1},  for  $x\in D_{t_0/4}(z_0)$ with $\rho_{z_0}(x)\leq t_0/(64L),$
$$\begin{aligned}
\P_x(\tau_D>t)&\leq c_{11}\P_x(X_{\tau_{U_1}}\in D)\\
&\leq c_{12}G_{U_1}(x, y_0)t_0^{d-2}\\
&\leq c_{12}G_{D_{t_0}(z_0)}(x, y_0)t_0^{d-2}\\
&\leq c_{12}G_{D_{t_0}(z_0)}(x, x_t)t_0^{d-2}
\leq c_{13}G_{D_{t_0}(z_0)}(x, x_t)t^{(d-2)/2},
\end{aligned}$$
where the fourth inequality is due to the standard chain argument and the  elliptic Harnack inequality \eqref{e:2.5}.
This establishes the upper bound of \eqref{e:4.19}.
The proof is complete.
\qed

\begin{thm}\label{T:5.1}
Suppose that $D$ is a Lipschitz open set in $\R^d$ with characteristics $(R_0, \Lambda_0).$
 Assume  the path distance in each connected component of $D$ is comparable to the Euclidean distance with characteristic $\chi_1.$
The followings are equivalent:

(i) For each $T>0,$ there exist positive constants $C_1=C_1(d, \phi, R_0, \Lambda_0, \chi_1, T)$ and $C_k=C_k(d, \phi,  \chi_1), k=2, 3$ such that for  any $x, y\in D$ and $t\in (0, T),$
\begin{equation}\label{e:5.7n}
C_1^{-1}\P_x(\tau_D>t)p(t, C_2x, C_2y)\P_y(\tau_D>t)\leq p_D(t, x, y)\leq C_1\P_x(\tau_D>t)p(t, C_3x, C_3y)\P_y(\tau_D>t).
\end{equation}

\smallskip
(ii) There exists $C_4=C_4(d, \phi, R_0, \Lambda_0)>1$ such that  for  any $z_0\in \partial D$ and $t\in (0, R_0^2),$
\begin{equation}\label{e:5.8n}
C_4^{-1} \E_x \tau_{D_{\sqrt t}(z_0)}\leq G_{D_{\sqrt t}(z_0)}(x, x_t)t^{d/2}\leq C_4 \E_x \tau_{D_{\sqrt t}(z_0)} \quad \hbox{for }\: x\in  D_{\sqrt t/4}(z_0),
\end{equation}
  where $x_t$ is a point of $D\cap \partial B(z_0, \sqrt t/2)$ with $\delta_D(x_t)\geq \kappa\sqrt t/2.$

\end{thm}

\pf  We first prove  \eqref{e:5.8n} $\Rightarrow $ \eqref{e:5.7n}. Let $T>0.$
It follows from Theorem \ref{T1} that there exists $c_0=c_0(d, \phi, T)>1$ such that $c_0^{-1} t^{-d/2}\leq p(t, x, y)\leq c_0 t^{-d/2}$ for $t\in (0, T)$ and $|x-y|\leq t^{1/2}.$
Thus by Theorem \ref{T2}, it only suffices to prove there exists $c=c(d, \phi, R_0, \Lambda_0, T)$ such that for $t\in (0, T),$ $x, y\in D$ with $|x-y|>\sqrt t,$
\begin{equation}\label{e:4.36}
\int_0^t \P_x(\tau_D>s)\P_y(\tau_D>t-s)\,ds\cdot j(|x-y|)
 \leq c\P_x(\tau_D>t)\P_y(\tau_D>t) p(t, x, y).
 \end{equation}

Let $t\in (0, T)$ and $t_0:=\sqrt t\wedge R_0.$ We will prove \eqref{e:4.36} in three separate cases.

(i) Suppose $|x-y|>\sqrt t$ with $\delta_D(x)<t_0/(8L)$ and $\delta_D(y)<t_0/(8L).$
It follows from \eqref{e:4.13n} that there exists $c_1=c_1(d, \phi, R_0, \Lambda_0, T)>0$ such that
\begin{equation}\label{e:4.24}\begin{aligned}
\int_0^t \P_x(\tau_D>s)\P_y(\tau_D>t-s)\,ds
\leq c_1\P_y(\tau_D>t)\E_x (\tau_D\wedge t)+c_1\P_x(\tau_D>t) \E_y (\tau_D\wedge t).
\end{aligned}\end{equation}
Let $z_x, z_y\in \partial D$ such that $\delta_D(x)=|x-z_x|$ and $\delta_D(y)=|y-z_y|.$
By Proposition \ref{P:2.2} (ii), there exists $c_2=c_2(d, \phi, R_0, \Lambda_0, T)>1$ such that for any  $t\in (0, T)$ and $s\in (0, t),$
$$\P_x(\tau_{U_{t_0}(z_x)}>s)\leq  \P_x(\tau_D>s) \leq c_2\P_x(\tau_{U_{t_0}(z_x)}>s).$$
Hence
\begin{equation}\label{e:5.7'}
\E_x (\tau_D\wedge t)=\int_0^t \P_x(\tau_D>s)\,ds\leq
c_2\int_0^t \P_x(\tau_{U_{t_0}(z_x)}>s)\,ds\leq c_2\E_x(\tau_{U_{t_0}(z_x)}\wedge t)\leq c_2\E_x \tau_{U_{t_0}(z_x)}.
\end{equation}
Similarly, \eqref{e:5.7'} holds for $y$ in place of $x.$
Then we have by \eqref{e:4.24} and \eqref{e:5.7'},
\begin{equation}\label{e:4.5}\begin{aligned}
\int_0^t \P_x(\tau_D>s)\P_y(\tau_D>t-s)\,ds
\leq c_1c_2\E_x\tau_{U_{t_0}(z_x)}\P_y(\tau_D>t)+c_1c_2\E_y\tau_{U_{t_0}(z_y)}\P_x(\tau_D>t).
\end{aligned}\end{equation}
Let $x_t$ be a point of $D\cap \partial B(z_x, \sqrt t_0/2)$ with $\delta_{D_{\sqrt t_0}(z_0)}(x_t)\geq \kappa\sqrt t_0/2.$
Suppose \eqref{e:5.8n} holds,  then there exists $c_3=c_3(d, \phi, R_0, \Lambda_0)$ such that $\E_x \tau_{D_{t_0}(z_x)}\leq  c_3 G_{D_{t_0}(z_x)}(x, x_t)t^{d/2}.$ Hence, by this inequality combining with Proposition \ref{L:4.5}, there exists $c_4=c_4(d, \phi, R_0, \Lambda_0, T)$ such that
\begin{equation}\label{e:4.6n}
 \E_x \tau_{U_{t_0}(z_x)}\leq  \E_x \tau_{D_{t_0}(z_x)}\leq  c_3 G_{D_{t_0}(z_x)}(x, x_t)t^{d/2}\leq c_4t\P_x(\tau_D>t).
  \end{equation}
Similarly, we have \eqref{e:4.6n} holds for $y$ in place of $x.$
In view of Theorem \ref{T1}, there exist $c_k=c_k(d, \phi, T)>0, k=5, 6$ such that
\begin{equation}\label{e:5.12}
 p(t, x, y)\geq c_5t^{-d/2}e^{-c_6|x-y|^2/t}+c_5tj(|x-y|) \geq c_5tj(|x-y|)  \quad \mbox{for} \quad |x-y|>\sqrt t.
\end{equation}
Hence by  \eqref{e:4.5}-\eqref{e:5.12}, for $t\in (0, T)$ and $|x-y|\geq \sqrt t,$
 $$\begin{aligned}
 &\int_0^t \P_x(\tau_D>s)\P_y(\tau_D>t-s)\,ds\cdot j(|x-y|) \\
  \leq &2c_1c_2c_4\P_x(\tau_D>t)\P_y(\tau_D>t)\cdot tj(|x-y|) \\
 \leq &2c_1c_2c_4c_5^{-1}\P_x(\tau_D>t)\P_y(\tau_D>t) p(t, x, y).
 \end{aligned}$$
 That is, \eqref{e:4.36} holds in this case.

 (ii) Suppose $|x-y|\geq \sqrt t$ with $\delta_D(x)\leq t_0/(8L)$ and $\delta_D(y)> t_0/(8L).$
   By Lemma \ref{L:2.7}, there exists $c_6=c_6(d, \phi, R_0, T)>0$ such that for $t\in (0, T),$
\begin{equation}\label{e:4.28}
1\geq \P_y(\tau_D>t)\geq \P_y(\tau_{B(y, t_0/(16L))}>t)\geq c_6.
\end{equation}
Note that $\E_y (\tau_D\wedge t)\leq t.$
Thus by \eqref{e:4.24}, \eqref{e:5.7'} and \eqref{e:4.6n},
\begin{equation}\label{e:4.29}
\begin{aligned}
\int_0^t \P_x(\tau_D>s)\P_y(\tau_D>t-s)\,ds
&\leq c_1\P_x(\tau_D>t) \E_y (\tau_D\wedge t)+c_1\P_y(\tau_D>t)\E_x (\tau_D\wedge t)\\
& \leq c_1\P_x(\tau_D>t)\cdot t+c_1c_2\E_x \tau_{U_{t_0}(z_x)}\\
& \leq c_1(1+c_2c_4)\P_x(\tau_D>t)\cdot t.
\end{aligned}\end{equation}
Then by \eqref{e:4.29} and \eqref{e:5.12}-\eqref{e:4.28},
$$\begin{aligned}
 &\int_0^t \P_x(\tau_D>s)\P_y(\tau_D>t-s)\,ds\cdot j(|x-y|) \\
 \leq & c_1(1+c_2c_4)\P_x(\tau_D>t)\cdot t j(|x-y|)\\
 \leq &c_1(1+c_2c_4)c_5^{-1}\P_x(\tau_D>t)\P_y(\tau_D>t) p(t, x, y).
 \end{aligned}$$
 Hence \eqref{e:4.36} holds for $|x-y|\geq \sqrt t$ with $\delta_D(x)\leq t_0/(8L)$ and $\delta_D(y)> t_0/(8L).$
 By the symmetry of $p_D(t, x, y)$ in $x$ and $y$,  \eqref{e:4.36} holds for  $|x-y|>\sqrt t$ with $\delta_D(x)> t_0/(8L)$ and $\delta_D(y)\leq t_0/(8L).$

(iii) Suppose  $|x-y|\geq \sqrt t$ with $\delta_D(x)> t_0/(8L)$ and $\delta_D(y)> t_0/(8L).$
 It follows from \eqref{e:5.12} and \eqref{e:4.28} that
 $$\begin{aligned}
 &\int_0^t \P_x(\tau_D>s)\P_y(\tau_D>t-s)\,ds\cdot j(|x-y|)\\
 &\leq t j(|x-y|)\leq c_5^{-1}p(t, x, y)\\
 &\leq c_5^{-1}c_6^{-2}\P_x(\tau_D>t)\P_y(\tau_D>t)p(t, x, y).
 \end{aligned}$$
 Hence \eqref{e:4.36} holds in this case.

\medskip

  Next we shall prove \eqref{e:5.7n} $\Rightarrow$ \eqref{e:5.8n}. Let $z_0\in \partial D$ and $t\in (0, R_0^2).$ Let $x\in D_{\sqrt t/4}(z_0)$ and $x_t$ be a point of $D\cap \partial B(z_0, \sqrt t/2)$ with $\delta_D(x_t)\geq \kappa \sqrt t/2.$
 Note that  $G_{D_{\sqrt t}(z_0)}(x, \cdot)$ is harmonic in $B(x_t, \kappa\sqrt t/4).$
By the elliptic Harnack inequality \eqref{e:2.5}, there exists $c_7=c_7(d, \phi)>0$ such that
 \begin{equation}\label{e:4.33}
 \E_x\tau_{D_{\sqrt t}(z_0)}\geq \int_{B(x_t, \kappa \sqrt t/8)} G_{D_{\sqrt t}(z_0)}(x, u)\,du\geq c_7G_{D_{\sqrt t}(z_0)}(x, x_t)t^{d/2}.
 \end{equation}
 Thus it suffices to prove that there exists $c_8=c_8(d, \phi, R_0, \Lambda_0)$ such that for $x\in D_{\sqrt t/4}(z_0),$
 \begin{equation}\label{e:4.44}
\E_x\tau_{D_{\sqrt t}(z_0)}\leq  c_8G_{D_{\sqrt t}(z_0)}(x, x_t)t^{d/2}.
 \end{equation}
 Let  $y_0$ be a point of $D\cap \partial B(z_0, R_0)$ with $\delta_D(y_0)>\kappa R_0.$
By Lemma \ref{L:2.7}, there exists $c_9=c_9(d, \phi, \kappa)>0$ such that for $s\in (0, R_0^2),$
  \begin{equation}\label{e:5.17}
  \P_{y_0}(\tau_D>s)\geq \P_{y_0}(\tau_D>R_0^2)\geq \P_{y_0}(\tau_{B(y_0, \kappa R_0)}>R_0^2)\geq c_9.
  \end{equation}
 Note that $|x-y_0|\geq R_0/2>\sqrt t/2$ for $x\in D_{\sqrt t/4}(z_0).$
 By an argument similar to that of Lemma \ref{L:4.4} with $(\sqrt t\wedge R_0)/2$ in place of $\sqrt t\wedge R_0$ there and \eqref{e:5.17}, for $x\in D_{\sqrt t/4}(z_0),$
 \begin{equation}\label{e:4.31}
 \begin{aligned}
 p_D(t, x, y_0)&\geq c_{10}\int_0^t \P_x(\tau_D>s)\P_{y_0}(\tau_D>t-s)\,ds\cdot j(|x-y_0|)\\
 &\geq c_9c_{10}\int_0^t \P_x(\tau_D>s)\,ds \cdot j(|x-y_0|)\\
 &=c_9c_{10}\E_x (\tau_D\wedge t) \cdot j(|x-y_0|)\\
 &\geq c_{11}\E_x \tau_{D_{\sqrt t}(z_0)} \cdot j(|x-y_0|)
 \end{aligned}  \end{equation}
 where the last inequality holds due to  $\E_x\tau_{D_{\sqrt t}(z_0)}\leq \E_x\tau_{B(x, 2\sqrt t)}\leq ct$ by Lemma \ref{L:2.2}.
 Since $|x-y_0|\geq R_0/2,$
 in view of Theorem \ref{T1} and Lemma \ref{L:4.3'},
  there exists $c_{12}=c_{12}(d, \phi)$ such that
 $p(t, x, y_0) \leq c_{12}tj(|x-y_0|)$ for any $t\in (0, 1).$
 Suppose  \eqref{e:5.7n} holds,  there exists $c_{13}=c_{13}(d, \phi, R_0, \Lambda_0)>0$ such that for $t\in (0, R_0^2)$ and $x\in D_{\sqrt t/4}(z_0),$
\begin{equation}\label{e:4.35}
p_D(t, x, y_0)\leq  c_{13}\P_x(\tau_D>t) p(t, x, y_0)\leq c_{13}c_{12}\P_x(\tau_D>t)tj(|x-y_0|).
\end{equation}
Thus by  \eqref{e:4.31} and \eqref{e:4.35},
 $$\dfrac{\E_x\tau_{D_{\sqrt t}(z_0)}}{t}\leq c_{11}^{-1}c_{13}c_{12}\P_x(\tau_D>t), \quad x\in D_{\sqrt t/4}(z_0).$$
 By Proposition \ref{L:4.5}, there exists $c_{14}=c_{14}(d, \phi, R_0, \Lambda_0)>0$ such that
$$ \P_x(\tau_D>t)\leq c_{14}G_{D_{\sqrt t}(z_0)}(x, x_t)t^{(d-2)/2}, \quad x\in D_{\sqrt t/4}(z_0).$$
The two inequalities above  deduces \eqref{e:4.44}.
This completes the proof.

\qed

\begin{thm}\label{T:5.2}
Suppose that $D$ is a Lipschitz open set in $\R^d$ with characteristics $(R_0, \Lambda_0).$
The followings are equivalent:

(i) There exists $C=C(d, \phi, R_0, \Lambda_0)>1$ such that  for  any $z_0\in \partial D$ and $t\in (0, R_0^2),$
\begin{equation}\label{e:5.23n}
C^{-1} \E_x \tau_{D_{\sqrt t}(z_0)}\leq G_{D_{\sqrt t}(z_0)}(x, x_t)t^{d/2}\leq C \E_x \tau_{D_{\sqrt t}(z_0)} \quad \hbox{for }\: x\in  D_{\sqrt t/4}(z_0),
\end{equation}
  where $x_t$ is a point of $D\cap \partial B(z_0, \sqrt t/2)$ with $\delta_D(x_t)\geq \kappa\sqrt t/2.$

\medskip

(ii) The scale invariant boundary Harnack principle holds for $X$ in  $D$.
That is,  there exists $ C= C(d, \phi, R_0,  \Lambda_0)>0$  such that for any $z_0\in\partial D, r\in (0, R_0/2)$ and any two  nonnegative harmonic functions $h_1$ and $h_2$ with respect to   $X$  in
 $D\cap B(z_0, 2r)$  vanishing continuously on $D^c\cap B(z_0, 2r),$
\begin{equation}\label{e:5.24n}
\dfrac{h_1(x)}{h_1(y)}\leq C \dfrac{h_2(x)}{h_2(y)} \quad \hbox{for any }
x,y\in D\cap B(z_0, r).
\end{equation}
\end{thm}

\pf  We first prove \eqref{e:5.23n} $\Rightarrow $ \eqref{e:5.24n}.
The proof is similar to \cite[Theorem 3.12]{CW2}.
 Let $z_0\in \partial D$ and $r\in (0, R_0/2).$ Let $h$ be a nonnegative harmonic function
on $D\cap B(z_0, 2r)$ and vanishes continuously on  $D^c\cap B(z_0, 2r).$
By the uniform Harnack inequality \eqref{e:2.5} and the standard chain argument,  it suffices to prove \eqref{e:5.24n} holds for $x\in D_{r/4}(z_0).$

Suppose \eqref{e:5.23n} holds. Let $x_0 \in D\cap  \partial B(z_0, r)$
with $\kappa r< \delta_D(x_0)< r.$ By a similar argument of  Lemma \ref{L:2.6} with $D_{2r}(z_0)$ in place of $U_r(z_0)$ there,
 there exists $c_1=c_1(d, \phi, \kappa)>0$  such that for any $z_0\in \partial D$ and   $r\in (0, R_0/2),$
\begin{equation}\label{e:4.37}
\P_x \left( X_{\tau_{D_{2r}(z_0)\setminus \overline{B(x_0, \kappa r/2)}}}\in \overline{B(x_0, \kappa r/2)} \right)
\geq c_1 r^{d-2} G_{D_{2r}(z_0)} (x, x_0)\geq c_2\dfrac{\E_x \tau_{D_{2r}(z_0)}}{r^2} \quad  \hbox{for } x\in D_{r/4}(z_0),
\end{equation}
where in the last inequality we used \eqref{e:5.23n}.
Let $\{U_n; n\geq 1\}$ be an increasing sequence of relatively compact subsets of $D_{2r}(z_0)$ so that
$\overline U_n\subset U_{n+1}$ and $\cup_{n=1}^\infty U_n = D_{2r}(z_0)$. Let $F:=D_{2r}(z_0)\setminus \overline{B(x_0, \kappa r/2)}.$
Then by \eqref{e:4.37}, for  $x\in D_{r/4 }(z_0)\subset F$,
\begin{equation}\label{e:5.22}\begin{aligned}
h(x)&= \lim_{n\to \infty} \E_x [ h(X_{\tau_{U_n\cap F}}) ]\\
&\geq   \lim_{n\to \infty}\E_x [h(X_{\tau_{U_n\cap F}}); X_{\tau_{U_n\cap F}}\in B(x_0, \kappa r/2)]\\
&\geq  \inf_{z\in \overline{B(x_0, \kappa r/2)}}h(z)\P_x (X_{\tau_{ F}}\in \overline{B(x_0, \kappa r/2)}).\\
 \end{aligned}\end{equation}
By \eqref{e:5.22}, \eqref{e:4.37} and  the uniform Harnack inequality \eqref{e:2.5}, there exists $c_3=c_3(d, \phi, R_0, \Lambda_0)>0$  such that
\begin{equation}\label{e:5.23}\begin{aligned}
h(x)
&\geq c_3h(x_0)\P_x \big(X_{\tau_{D_{2r}(z_0)\setminus \overline{B(x_0, \kappa r/2)}}}\in \overline{B(x_0, \kappa r/2)} \big)\\
&\geq c_2c_3h(x_0) r^{-2}  \E_x  \left[ \tau_{D_{r/2}(z_0)} \right].
\end{aligned}\end{equation}

On the other hand, by Lemma \ref{L:2.1'}, for $x\in D_{r/4}(z_0),$
\begin{equation}\label{e:5.29'}\begin{aligned}
 h(x)
&=\E_x h(X_{\tau_{D_{r/2}(z_0)}})\\
&=\E_x  \left[h(X_{\tau_{D_{r/2}(z_0)}}); X_{\tau_{D_{r/2}(z_0)}}\in D_{r}(z_0) \right] \\
&\quad +\E_x \left[ h(X_{\tau_{D_{r/2}(z_0)}}); X_{\tau_{D_{r/2}(z_0)}}\in  B^c(z_0, r) \right].
\end{aligned}\end{equation}
 By the Carleson estimate Proposition \ref{L:2.4}  and  \cite[Lemma 2.5]{CW2},  there exists $c_4=c_4(d, \phi, R_0, \Lambda_0)$ such that
 \begin{eqnarray}\label{e:5.30'}
   \E_x [h(X_{\tau_{D_{r/2}(z_0)}}); X_{\tau_{D_{r/2}(z_0)}}\in D_{r}(z_0)]
 &\leq &\sup_{z\in D_{r}(z_0)}h(z)\cdot\P_x(X_{\tau_{D_{r/2}(z_0)}}\in  D_{r}(z_0)) \nonumber \\
&\leq &c_4h(x_0) r^{-2}  \E_x  \left[ \tau_{D_{r/2}(z_0)}\right]  .
 \end{eqnarray}
For the second item in \eqref{e:5.29'}, by the condition \eqref{e:1.1}, there exists
 $c_5=c_5(d, \phi)$
such that
\begin{eqnarray*}
&&\E_x  \big[ h(X_{\tau_{D_{r/2}(z_0)}}); X_{\tau_{D_{r/2}(z_0)}}\in B^c(z_0, r) \big]\\
&=& \E_x\int_0^{\tau_{D_{r/2}(z_0)}}\int_{B^c(z_0, r)} h(z) j(|X_s-z|) \,dz\,ds\\
&\leq & c_5  \E_x \left[ \tau_{D_{r/2}(z_0)} \right] \int_{ B^c(z_0, r)}j(|z-z_0|) h(z)\,dz.
\end{eqnarray*}
Let $y_0 \in D\cap \partial B(z_0, r/2)$ with $\kappa r/2< \delta_D(y_0)< r/2.$
It follows from Lemma \ref{L:2.2} and the condition \eqref{e:1.1} that
\begin{eqnarray*}
h(y_0) &\geq &  \E_{y_0} \big[ h(X_{\tau_{B(y_0, \kappa r/2)}}); X_{\tau_{B(y_0, \kappa r/2)}}\in B^c(z_0, r) \big]\\
&= &\E_{y_0}\int_0^{\tau_{B(y_0, \kappa r/2)}}\int_{B^c(z_0, r)} h(z) j(|X_s-z|) \,dz\,ds\\
&\geq&  c_6\E_{y_0}\left[ \tau_{B(y_0, \kappa r/2)} \right]  \int_{ B^c(z_0, r)}j(|z-z_0|) h(z)\,dz\\
&\geq & c_7 r^2   \int_{ B^c_{r}(z_0)}j(|z-z_0|) h(z)\,dz.
\end{eqnarray*}
By the standard Harnack chain and the elliptic Harnack principle  \eqref{e:2.5}, there exists $c_8=c_8(d, \phi, \Lambda_0)>0$ such that $h(y_0)\leq c_8 h(x_0).$
Hence,
\begin{equation}\label{e:5.31}
\E_x  \left[ h(X_{\tau_{D_{r/2}(z_0)}}); X_{\tau_{D_{r/2}(z_0)}}\in B^c(z_0, r) \right]
 \leq  \dfrac{c_5c_8}{c_7}h(x_0) r^{-2}  \E_x \left[ \tau_{D_{r/2}(z_0)}\right].
\end{equation}
Hence by combining \eqref{e:5.29'}-\eqref{e:5.31},
\begin{equation}\label{e:5.24}
h(x)\leq c_9h(x_0)r^{-2} \E_x  \big[\tau_{D_{r/2}(z_0)} \big],
\end{equation}
where $c_9=c_4+\frac{c_5c_8}{c_7}.$
Thus by combining \eqref{e:5.23} and \eqref{e:5.24}, there exists $c_{10}=c_{10}(d, \phi, R_0, \Lambda_0)>1$ such that
\begin{equation}\label{e:5.29}
c_{10}^{-1} \dfrac{\E_x  \big[ \tau_{D_{r/2}(z_0)} \big] }{r^2}\leq \dfrac{h(x)}{h(x_0)}\leq c_{10} \dfrac{\E_x  \big[ \tau_{D_{r/2}(z_0)} \big] }{r^2} \quad \mbox{for} \quad x\in  D_{r/4}(z_0).
\end{equation}
This yields that \eqref{e:5.24n} holds for $x\in D_{r/4}(z_0)$.

\smallskip
Next, we  prove   \eqref{e:5.24n} $\Rightarrow $ \eqref{e:5.23n}.
Suppose the scale invariant boundary Harnack principle \eqref{e:5.24n} holds in a Lipschitz open set $D$.
Note that the transition density function of $X$ in $\R^d$ exists and the condition \eqref{e:1.1} holds, thus the assumptions (A1) and (A3) in \cite{CW2} hold.
Let $z_0\in \partial D$ and $t\in (0, R_0^2).$
 Let $y_0$ be a point of $D\cap \partial B(z_0, \sqrt t/4)$ with
$\delta_{D_{\sqrt t}(z_0)}(y_0)\geq \kappa\sqrt t/4.$
By Theorem 1.3 (ii) in \cite{CW2}, if the scale invariant boundary Harnack principle  \eqref{e:5.24n}  holds, then there exists $c_{11}=c_{11}(d, \phi, R_0, \Lambda_0)>0$  such that
for $t\in (0, R_0^2),$
\begin{equation}\label{e:5.30}
\dfrac{\E_x \tau_{D_{\sqrt t}(z_0)}}{\E_{y_0} \tau_{D_{\sqrt t}(z_0)}}
\leq c_{11}\dfrac{G_{D_{\sqrt t}(z_0)}(x, x_t)}{G_{D_{\sqrt t}(z_0)}(y_0, x_t)}, \quad x\in  D_{\sqrt t/4}(z_0),
\end{equation}
 where $x_t$ is a point of $D\cap \partial B(z_0, 3\sqrt t/4)$ with $\delta_D(x_t)\geq \kappa\sqrt t/2.$
By Lemma \ref{L:2.2}, there exists $c_{12}=c_{12}(d, \phi)>0$ such that $\E_{y_0} \tau_{D_{\sqrt t}(z_0)}\leq \E_{y_0} \tau_{B(y_0, 2\sqrt t)}\leq c_{12}t.$
Moreover, by the Harnack chain argument and Lemma \ref{L:2.4'}, there exists $c_k=c_k(d, \phi, \Lambda_0)>0, k=13, 14$ such that $G_{D_{\sqrt t}(z_0)}(y_0, x_t)\geq c_{13}\inf_{y\in B(x_t, \kappa \sqrt t/4)}
G_{B(x_t, \kappa \sqrt t/2)}(y, x_t)\geq c_{14} t^{(2-d)/2}.$
By applying these two inequalities to \eqref{e:5.30}, the inequality on the left side of \eqref{e:5.23n} holds.
The inequality on the right side of \eqref{e:5.23n} follows from \eqref{e:4.33}.
\qed

\smallskip

 {\bf Proof of Theorem \ref{T3}.} Theorem \ref{T3} is obtained by Theorems \ref{T:5.1} and \ref{T:5.2}.

 \qed

\bigskip

In the last part of this paper, we shall prove Corollary \ref{C2}.
Recall that an open set $D$ in $\R^d$ (when $d\geq 2$) is said to be $C^{1,1}$ if there exist a localization radius
$R_0 >0$ and a constant $\Lambda_0>0$ such that for every $z\in\partial D,$ there exist a $C^{1, 1}$ function
$\Gamma=\Gamma_z: \R^{d-1}\rightarrow \R$ satisfying
\begin{equation}\label{e:1.5}
\Gamma(0)=\nabla \Gamma(0)=0, \quad \|\nabla \Gamma\|_\infty\leq \Lambda_0, \quad
|\nabla \Gamma(x)-\nabla \Gamma(y)|\leq \Lambda_0 |x-y|
\end{equation}
 and an orthonormal coordinate system
$CS_z: y=(y^1, \cdots, y^{d-1}, y^d)=:(\wt {y}, y^d)\in \R^{d-1}\times \R$ with its origin at $z$ such that
$$
B(z, R_0)\cap D=\{y=(\wt {y}, y^d)\in B(0,R_0) \hbox{ in } CS_z: y^d>\Gamma(\wt {y})\}.
$$
The pair $(R_0, \Lambda_0)$ is called the characteristics of the $C^{1, 1}$ open set $D.$

\medskip

 {\bf Proof of Corollary \ref{C2}.} Let $d\geq 3$ and $D$ be a $C^{1, 1}$ open set  with characteristics $(R_0, \Lambda_0).$
By  \cite[Lemma 2.2]{S}, for each $r\in (0, R_0),$ there exists $L=L(R_0, \Lambda_0, d)>0$
such that for any $z\in\partial D,$ there is a $C^{1,1}$ connected open set  $U_{z, r}\subset D$  with characteristics $(rR_0/L, L\Lambda_0/r)$ such that $D\cap B(z, r)\subset U_{z, r}\subset D\cap B(z, 2r).$

It follows from \cite[(2.19) and Lemma 2.17]{CW} that there exists $c_1=c_1(d, \phi, R_0, \Lambda_0)>1$ such that for any $z\in \partial D, r\in (0, R_0)$ and $x, y\in U_{z, r}$ with $y\neq x,$
\begin{equation}\label{e:5.6}
\begin{aligned}
&c^{-1}_1|x-y|^{2-d} \left(1\wedge \frac{\delta_{U_{z, r}}(x)}{|x-y|}\right) \left(1\wedge \frac{\delta_{U_{z, r}}(y)}{|x-y|}\right)\\
&\quad \leq G_{U_{z, r}}(x, y)\leq c_1|x-y|^{2-d} \left(1\wedge \frac{\delta_{U_{z, r}}(x)}{|x-y|}\right) \left(1\wedge \frac{\delta_{U_{z, r}}(y)}{|x-y|}\right).
\end{aligned}\end{equation}
 Note that $D$ is also a Lipschitz open set  with characteristics $(R_0, \Lambda_0).$ Recall that $\kappa=\kappa(d, \Lambda_0)$ is the constant in \eqref{e:1.8}. For each $z\in \partial D,$ let $z_{r/2}$ be a point of $D\cap \partial B(z, r/2)$ with $\delta_D(z_{r/2})\geq \kappa r/2.$
 By \eqref{e:5.6}, there exists $c_2=c_2(d, \phi, R_0, \Lambda_0)>1$ such that for $z\in \partial D$ and $r\in (0, R_0).$
\begin{equation}\label{e:6.37}
c_2^{-1} r^{1-d}\delta_D(x)\leq G_{U_{z, r}}(x, z_{r/2}) \leq c_2 r^{1-d}\delta_D(x) \quad  \mbox{for} \quad x\in D_{r/4}(z),
\end{equation}
By \eqref{e:5.6} and $\E_x \tau_{U_{z, r}}=\int_{U_{z, r}}G_{U_{z, r}}(x, y)\,dy$, there exists $c_3=c_3(d, \phi, R_0, \Lambda_0)>1$ such that for $z\in \partial D$ and $r\in (0, R_0),$
$$c_3^{-1} \delta_D(x)r\leq \E_x \tau_{U_{z, r}}\leq c_3 \delta_D(x)r \quad  \mbox{for} \quad x\in D_{r/2}(z).$$
Note that  $U_{z, r/2}\subset D_r(z)\subset U_{z, r}.$
By \eqref{e:6.37} and the elliptic Harnack inequality \eqref{e:2.5}, there exists $c_4=c_4(d, \phi, R_0, \Lambda_0)>1$ such that
\begin{equation}\label{e:6.40}
c_4^{-1}r^{1-d}\delta_D(x)\leq G_{D_r(z)}(x, z_{r/2})\leq c_4r^{1-d}\delta_D(x) \quad  \mbox{for} \quad x\in D_{r/4}(z).
\end{equation}
It is easy to see that for $z\in \partial D$ and $r\in (0, R_0),$
\begin{equation}\label{e:6.38'}
(2c_3)^{-1} \delta_D(x)r\leq \E_x \tau_{D_r(z)}\leq c_3\delta_D(x) r \quad  \mbox{for} \quad x\in D_{r/4}(z).
\end{equation}
Thus by \eqref{e:6.38'} and \eqref{e:6.40}, for $z\in \partial D$ and $r\in (0, R_0),$
 \begin{equation}\label{e:5.7}
 G_{D_r(z)}(x, z_{r/2})r^d \asymp \E_x \tau_{D_r(z)}, \quad x\in  D_{r/4}(z),
 \end{equation}
 where the comparison constant depends on $(d, \phi, R_0, \Lambda_0).$
 Then by \eqref{e:5.7} and Theorem \ref{T3}, there exists $c_k, k=5, 6, 7$  such that for any $t\in (0, T)$ and $x, y\in D,$
 \begin{equation}\label{e:5.8}
c_5^{-1}\P_x(\tau_D>t)\P_y(\tau_D>t)p(t, c_6x, c_6y)\leq  p_D(t, x, y)\leq c_5\P_x(\tau_D>t)\P_y(\tau_D>t)p(t, c_7x, c_7y).
 \end{equation}
By Proposition \ref{L:4.5} and \eqref{e:6.40} with $r=\sqrt t\wedge R_0$, there exists $c_8=c_8(d, \phi, R_0,  \Lambda_0, T)>1$ such that for any $z_0\in \partial D$ and $t\in (0, T),$
 \begin{equation}\label{e:5.9}
 c_8^{-1}\dfrac{\delta_D(x)}{\sqrt t}\leq \P_x(\tau_D>t)\leq c_8\dfrac{\delta_D(x)}{\sqrt t}, \quad x\in  D_{(\sqrt t\wedge R_0)/4}(z_0).
 \end{equation}
  If $x\in D$ with $\delta_D(x)> (\sqrt t\wedge R_0)/4,$ then by Lemma \ref{L:2.7}, there exists $c_9=c_9(d, \phi, R_0,  T)\in (0, 1)$ such that for $t\in (0, T),$
\begin{equation}\label{e:5.10}
1\geq \P_x(\tau_D>t)\geq \P_x(\tau_{B(x, (\sqrt t\wedge R_0)/8)}>t)\geq c_9.
\end{equation}
Hence, by \eqref{e:5.9} and \eqref{e:5.10},
there exists $c_{10}=c_{10}(d, \phi, R_0, \Lambda_0, T)>1$ such that for any $t\in (0, T),$
 \begin{equation}\label{e:6.9}
 c_{10}^{-1}\left(1\wedge \dfrac{\delta_D(x)}{\sqrt t}\right)\leq \P_x(\tau_D>t)\leq c_{10}\left(1\wedge \dfrac{\delta_D(x)}{\sqrt t}\right), \quad x\in  D.
 \end{equation}
 By combining \eqref{e:5.8}, \eqref{e:6.9} and Theorem \ref{T1},  \eqref{e:1.14} can be obtained.
When $D$ is a bounded $C^{1, 1}$ domain, for large time $t\geq 3$, \eqref{e:1.15} is obtained by \eqref{e:6.9} and Theorem \ref{T2}(iii).

\qed

\medskip

{\bf Acknowledgement.} The author is  grateful to Professor Zhen-Qing Chen for the valuable and helpful comments.

\end{document}